\documentclass[final]{siamart1116}

\usepackage{color}
\usepackage{epsfig,bm,amsfonts}
\usepackage{graphicx, caption}
\usepackage{subcaption}
\usepackage{epstopdf, tikz}
\usepackage{algorithmic,amsopn}
\usepackage{booktabs}
\usepackage{arydshln}
\usepackage{url}
\usetikzlibrary{shapes,arrows}
\usetikzlibrary{arrows,decorations.pathmorphing,backgrounds,fit,positioning,shapes.symbols,chains}

\DeclareMathAlphabet\mathbold{OML}{cmm}{b}{it}
\numberwithin{equation}{section}

\newtheorem*{remark}{Remark}

\newtheorem{problem}{Problem}[section]

\renewcommand{\div}{\operatorname{div}}

\newcommand{\MFEM}{{\sc MFEM}}
\newcommand{\MOONOLITH}{{\sc MOONoLith}}

\def\bu{{\mathbf u}}
\def\bv{{\mathbf v}}
\def\bp{{\mathbf p}}
\def\bq{{\mathbf q}}
\def\bx{{\mathbf x}}
\def\by{{\mathbf y}}
\def\bn{{\mathbf n}}
\def\f{{\mathbf f}}
\def\bX{{\mathbf X}}
\def\bR{{\mathbf R}}
\def\oq{{\overline q}}
\def\obu{{\overline \bu}}
\def\obv{{\overline \bv}}
\def\otheta{{\overline \theta}}
\def\obtheta{\overline{\boldsymbol \theta}}
\def\btheta{{\boldsymbol \theta}}
\def\oD{{\overline D}}
\def\blambda{{\boldsymbol \lambda}}

\newcommand{\Mat}{Mat\'{e}rn }
\newcommand{\V}{\mathbb{V}}
\newcommand{\E}{\mathbb{E}}
\newcommand{\KL}{Karhunen-Lo\`{e}ve\ }
\newcommand{\norm}[1]{\left\Vert#1\right\Vert}
\newcommand{\abs}[1]{\left\vert#1\right\vert}
\newcommand{\curly}[1]{\left\{#1\right\}}
\newcommand{\p}[1]{\left(#1\right)}
\newcommand{\brac}[1]{\left[#1\right]}

\newcommand{\oP}{{\overline P}}
\newcommand{\R}{{\mathbb R}}
\newcommand{\oR}{{\overline {\mathbf R}}}
\newcommand{\oTheta}{{\overline \Theta}}
\newcommand{\T}{{\mathcal T}}
\newcommand{\oT}{{\overline \T}}
\newcommand{\oF}{{\overline F}}
\newcommand{\oh}{{\bar h}}
\newcommand{\otau}{{\overline \tau}}
\newcommand{\ovarphi}{{\overline \varphi}}
\newcommand{\os}{{\overline s}}

\newcommand{\A}{{\mathcal A}}

\newcommand{\M}{{\mathcal M}}
\newcommand{\W}{\mathcal{W}}
\renewcommand{\P}{{\mathcal P}}

\DeclareMathOperator{\diag}{diag}

\title{Scalable hierarchical PDE sampler for generating spatially correlated
    random fields using non-matching meshes
    \footnote{ S\lowercase{ubmitted: 19-}M\lowercase{ay-2017} A\lowercase{ccepted: 5-}D\lowercase{ec-2017} }
    \footnote{T\lowercase{his work is performed under the auspices of the} U.S.
    D\lowercase{epartment of} E\lowercase{nergy under} C\lowercase{ontract}
    DE-AC52-07NA27344. LLNL-JRNL-731006.}}

\author{Sarah Osborn \thanks{Center for Applied Scientific Computing,
    Lawrence Livermore National Laboratory, P.O. Box 808, L-561, Livermore,
    CA 94551, USA.
    (osborn9@llnl.gov).}
    \and Patrick Zulian\thanks{Institute of Computational Science, Universit\`{a} della Svizzera
    italiana, 6900, Lugano, Switzerland.}
    \and Thomas Benson \footnotemark[3]
    \and Umberto Villa \thanks{Institute for Computational Engineering and
    Sciences, University of Texas, Austin, TX.}
    \and Rolf Krause\footnotemark[4]
    \and Panayot S. Vassilevski \footnotemark[3] \thanks{Fariborz Maseeh Department of Mathematics and Statistics, Portland State University, Portland, OR.}}

\begin{document}
\maketitle
\begin{abstract}
This work describes a domain embedding technique between two non-matching meshes
used for generating realizations of spatially correlated random fields with applications to
large-scale sampling-based uncertainty quantification.
The goal is to apply the multilevel Monte Carlo (MLMC) method
for the quantification of output uncertainties of PDEs with random input
coefficients
on general, unstructured computational domains.
We propose a highly scalable, hierarchical sampling method to generate
realizations of a Gaussian random field on a given unstructured mesh by solving
a reaction-diffusion PDE with a stochastic right-hand side.
The stochastic PDE is discretized using the mixed
finite element method on an embedded domain with a structured
mesh, and then the solution is projected onto the unstructured mesh.
This work describes implementation details on how to efficiently transfer data
from the structured and unstructured meshes at coarse levels, assuming this can
be done efficiently on the finest level. 
We investigate the efficiency and parallel scalability of the technique for the
scalable generation of Gaussian random fields in three dimensions. 
An application of the MLMC method is presented for quantifying uncertainties of
subsurface flow problems. We demonstrate the scalability of the sampling method with non-matching 
mesh embedding, coupled with a parallel forward model problem solver, for large-scale 3D 
MLMC simulations with up to $1.9\cdot 10^9$ unknowns.
\end{abstract}


\section{Introduction} \label{sec:intro} 
Many mathematical models of physical phenomena involve spatially varying input
data which is often subject to uncertainty. This uncertainty will propagate
through a simulation and lead to uncertainty in the output. The goal in
forward propagation of uncertainty is to quantify the
effect of the input uncertainties in the output of numerical simulations. 
We consider models based on
partial differential equations (PDEs) with spatially correlated input coefficients
subject to uncertainty that are modeled as a random field with particular
statistical properties. Then our goal is to compute statistics of the solution
to the PDE with random input coefficients for large-scale problems
using Monte Carlo methods. In particular, we consider the multilevel Monte Carlo
(MLMC)
method, which runs repeated simulations at random realizations of the uncertain
input data on a hierarchy of spatial resolutions. Then, the approximations are
used to compute corresponding sample averages of the desired statistics of the
solution of the PDE. 

As an example model problem, we consider the simulation of subsurface flows
governed by Darcy's law. 
The permeability tensor, $k$, is often subject to uncertainty, 
due to a lack of knowledge of the porous medium at all locations. 
To account for this uncertainty, the permeability field is 
modeled as a random field with given mean and covariance structure.
Estimating the impact of uncertainty on the results of a groundwater flow
simulation    
is useful in many situations, for example in risk analysis for radioactive waste disposal or in oil
reservoir simulations. 

Monte Carlo techniques are a widely used class of methods 
to estimate particular quantities of
interest for PDEs with random input coefficients, and require the solution of
the model problem equations for many different realizations of $k$.
The computational cost for large-scale problems can often be
prohibitively large, as computing each sample amounts to solving a PDE with a
fine mesh. 
MLMC methods~\cite{heinrich2001multilevel,giles2008multilevel} are 
used to accelerate the convergence of standard Monte Carlo methods and offer 
significant computational savings.
These methods employ a hierarchy
of spatial resolutions as a variance reduction technique for the approximation
of expected quantities of interest and have been successfully applied to a wide
variety of
applications; see, e.g.,~\cite{barth2011multi,
cliffe2011multilevel, teckentrup2013further, doi:10.1137/110853054,
Gittelson2013}. 
As in \cite{o16}, 
for our MLMC simulations we consider a general unstructured fine grid and construct a
hierarchy of algebraically coarsened grids and finite element spaces using 
element-based algebraic multigrid techniques (AMGe), which 
possess the same order approximation property as the original fine level
discretization, see e.g.,
\cite{lashuk2014construction,lashuk2012element,pasciak2008exact}.
Of particular importance is the ability to run MLMC simulations on large-scale 
unstructured meshes of complicated computational domains. 
To accomplish this task, input realizations of the random field must 
first be generated for general unstructured meshes. 
Thus, an important task for large-scale MLMC simulations is the 
scalable generation of Gaussian random field realizations, which is the focus of
this work.

A common choice in stochastic modeling for subsurface hydrology 
is to model the random permeability, 
$k$, as a log-normal random field, so that
$k(\bx,\omega) = \exp[\theta(\bx,\omega)]$, where $\theta(\bx,\omega)$ is a Gaussian
random field with prescribed mean and covariance structure; see,
e.g.,~\cite{delhomme1979spatial,gelhar1993}. 
Several methods exist to realize samples of Gaussian random fields to be used in
MLMC simulations.
The widely-used \KL (KL) expansion~\cite{loeve1978probability} provides an
infinite series representation of the random field involving the eigenvalues and
eigenvectors of the integral operator associated with the covariance function.
In practical computations, the series is truncated, which results in a
truncation error (bias) in
a Monte Carlo simulation. Additionally, the computation quickly becomes
infeasible
for large-scale simulations as a dense eigenvalue problem must be solved.
Approaches
based on randomized methods and hierarchical semi-separable matrices
can drastically reduce the cost of solving the eigenvalue problem
(see, e.g., \cite{SaibabaLeeKitanidis16}); however, only the
dominant eigenmodes of the KL expansion are computed and, therefore, introduce
bias in the sampling.
Circulant embedding \cite{dietrich1997, park2015block} 
offers a fast and exact 
simulation of stationary Gaussian random fields on a regular grid.
This method exploits the Fast Fourier Transform (FFT) method to implicitly
construct a basis that (block)-diagonalizes the covariance matrix.
However, the computational cost of the method depends 
on the correlation length of the random field. Additionally,
the random field is assumed to be stationary, whereas our proposed
method can handle correlation functions and marginal variance 
functions that are space-dependent. Last but not least, although
scalable implementations of three-dimensional FFT are available (see, e.g.
\cite{pekurovsky2012p3dfft, pippig2013pfft}), we are not aware of any
parallel publicly available implementations of circulant embedding for large-scale computations.
For these reasons, we do not pursue this approach.

An alternative technique for generating realizations of $\theta(\bx,\omega)$
relies on the link between Gaussian fields and Gaussian
Markov random fields,
where a stochastic
partial differential equation (SPDE) with a white noise forcing term is solved to generate the desired
realizations~\cite{whittle1954stationary,whittle1963stochastic,lindgren2011explicit}.
This approach provides a sampling method that is highly
scalable as the method leverages solution strategies for solving sparse linear
systems arising from the
finite element discretization of the SPDE, as investigated in \cite{o16}. 
A limitation of the SPDE sampling method is that the computed realizations contain artificial
boundary effects, arising from the discretization of the SPDE on a finite
domain. A possible solution is to embed the computational domain into a larger
one as investigated in \cite{o16}; however, this poses some challenges for unstructured meshes.

As a follow-up and as an alternative to the method proposed in~\cite{o16}, 
we propose a scalable domain embedding technique using non-matching meshes. 
The SPDE with white noise forcing term is discretized and solved on a regular, structured mesh,
then is projected back to the original, unstructured mesh of interest. 
We use a completely parallel approach that allows for the transfer of discrete
fields between unstructured volume and surface meshes, which can be arbitrarily
distributed among different processors~\cite{krause2016parallel}.
Then the resulting realization of the random field can be used as the input realization
of a Monte Carlo method. 

The key contribution of this work 
is to provide a flexible, black-box
workflow for embedding complex 3D domains in parallel 
for a highly scalable, hierarchical sampler of Gaussian random fields.
The domain embedding is necessary for the alleviation of
boundary artifacts in the SPDE sampler, and allows for the use of more efficient 
solvers for structured grids. 
In particular, 
we will use a scalable hybridization multigrid preconditioning strategy.
This sampling technique allows for the use of unstructured meshes for complex
computational geometries in 3D, which is necessary for realistic subsurface flow
simulations. 

The focus of this work is on making MLMC simulations feasible in practice
for large-scale problems. To this end, we focus on parallelism across the spatial domain in
computing realizations of the input random field using a novel technique, 
then performing the subsequent solve of the model of interest in parallel. 
Using the presented sampling strategy, coupled with scalable techniques for the solution
of the forward model problem, we demonstrate that the approach allows for the solution
of an extreme-scale
forward UQ problem with $1.9\cdot 10^9$ unknowns with high accuracy. 
Moreover, scalability of our approach can be further improved by exploiting additional
levels of parallelism, such as the scheduling approach
\cite{GmeinerDRSW16}, where the authors investigate the complex task
of scheduling parallel tasks within and across levels of MLMC.

The remainder of the paper is organized as follows. 
The standard Monte Carlo and MLMC methods
are reviewed in Section~\ref{sec:mlmc}. 
The forward model problem and discretization is described in Section~\ref{sec:darcy}.
In Section~\ref{sec:gaussian-rf} we
discuss a method for generating realizations of spatially correlated
random fields based on Gaussian Markov random fields. The sampling method is based on
solving a mixed discretization of a reaction-diffusion equation
with a stochastic right-hand side using domain embedding with two non-matching meshes. 
The scalable mapping of discrete fields
between non-matching meshes is discussed, and a brief overview of the 
implementation details of the projection operator is provided in Section~\ref{sec:projection}. 
The hierarchical SPDE sampling procedure is introduced and examined in 
Section~\ref{sec:ml-sampler}. Additionally, we discuss the iterative solution strategy solving for
the resulting saddle-point mixed systems.
Numerical results are presented in Section~\ref{sec:results} that investigate the
parallel performance of the proposed sampling
method, and an adaptive MLMC simulation for uncertainty quantification in
subsurface flow using different geometries in three spatial
dimensions.
Concluding remarks are given in Section~\ref{sec:conclusions}.

\section{Multilevel Monte Carlo Methods} \label{sec:mlmc}
In this section we briefly review the standard Monte Carlo (MC) method and the 
MLMC method for computing moments of quantities of interest 
 $Q(\omega)=\mathcal{B}\brac{\mathbf{X}(\bx,\omega)}$, where $\mathbf{X}(\bx,\omega)$
 is the solution of a PDE with random input coefficient following the
presentation in~\cite{cliffe2011multilevel}.
In our model problem, the quantity of interest is 
related to the pressure and/or Darcy flux of the mixed Darcy equations \eqref{eq:darcy}. 
In practice, the inaccessible quantity of interest $Q(\omega)$ is approximated by
$Q_{h}(\omega)$, the functional of the finite element solution
$\mathbf{X_h}(\bx,\omega)$ on the triangulation $\T_h$. 

\subsection{Standard Monte Carlo simulation} 
The standard Monte Carlo estimator
for $\E[Q]$ is  
\begin{equation} 
    \widehat{Q}_{h}^{MC} = \frac{1}{N}\sum_{i=1}^{N}Q_{h}^i, 
\end{equation} 
where $Q_{h}^i$ is the $i^{\text{th}}$ sample of $Q_{h}$ and $N$ is the number of
(independent) samples.

The mean square error (MSE) of the method is given by
\begin{equation}
    \E\brac{(\widehat{Q}_h^{MC}-\E[Q])^2} = \frac{1}{N}\V[Q_{h}]  + \p{\E[Q_h
    -Q]}^{2}.  
\end{equation}
Thus, the error naturally splits into two terms: the sampling error given 
by the variance of the estimator, and the estimator bias related 
to the finite element discretization error.
The estimator variance decays linearly with respect to the sample size $N$, 
and the bias gets smaller as the discretization is refined. This 
can make the method prohibitively expensive 
as the necessary samples size becomes very large and a fine
spatial discretization is necessary for high accuracy.

\subsection{Multilevel Monte Carlo simulation}
The MLMC method~\cite{cliffe2011multilevel,giles2008multilevel} 
is an effective variance reduction technique, which reduces the overall 
computational cost of the standard MC method using a hierarchical
sampling technique. 
Assume we have a sequence $Q_L,\dots,Q_1$ which approximates the 
quantity of interest $Q_0 = Q_h$ with
increasing accuracy and increasing cost. 
The sequence of approximations is often found by solving the
model problem on a geometric hierarchy of meshes constructed by uniform
refinement, but other alternative options have been considered; see,
e.g.,~\cite{Haji-Ali2016}. 
As in~\cite{o16}, we consider a nested hierarchy of spatial approximations 
constructed from AMGe methods for finite element discretizations, 
which possess the same order approximation property as the original fine level
discretization, discussed in Section~\ref{subsec:ml-darcy}. 

Using the linearity of the expectation operator, 
we have the following expression for $\E[Q_{h}]$ based on corrections with respect
to the next coarser discretization level:
\begin{equation} \E[Q_{h}] = \E[Q_{L}] +
    \sum_{\ell=0}^{L-1}\E[Q_{\ell}-Q_{\ell+1}] = \sum_{\ell=0}^{L}
    \E[Y_{\ell}], 
    \label{eq:telescope}
\end{equation} where $Y_{\ell}:=Q_{\ell}-Q_{\ell+1}$ for $i=0,\ldots L-1$ and
$Y_{L}:=Q_{L}$.
A standard MC estimator is used to independently estimate 
the expectation of $Y_{\ell}$ on each level, with suitably chosen samples sizes 
to minimize the overall computational complexity, yielding
the MLMC estimator for \eqref{eq:telescope} given by
\begin{equation}
    \widehat{Q}_{h}^{MLMC}
    =\sum_{\ell=0}^{L}\brac{\frac{1}{N_{\ell}}\sum_{i=1}^{N_{\ell}} Y_{\ell}^{(i)} }.
    \label{eq:mlmc_estimator}
\end{equation}
It is important to note that in \eqref{eq:mlmc_estimator} for a particular level $\ell$, the same
random sample $\omega^{(i)}$ is used with two spatial discretizations 
to compute $\bX_{\ell}(\bx,\omega^i)$ and
$\bX_{\ell+1}(\bx,\omega^i)$ when estimating the quantity $Y_{\ell}^{(i)}$.

The mean square error for the MLMC method becomes 
\begin{equation}
    \E\brac{\p{\widehat{Q}_{h}^{MLMC}-\E[Q]}^2}=
    \sum_{\ell=0}^{L}\frac{1}{N_{\ell}}\V[Y_{\ell}] +
    \p{\E\brac{Q_{0}-Q}}^2.  \label{eq:mlmc_mse} 
\end{equation} 
Similar to the standard MC error, the two terms of the MLMC MSE
represent the variance of the estimator and the discretization error. 

For a prescribed MSE of less than $\varepsilon^2$, 
the spatial discretization of the finest
level of the hierarchy is chosen so the bias term is less than $\varepsilon^2/2$.
Then the number of samples at each level $\ell$ is chosen to minimize the overall computational
cost leading to the following formula for the optimal number of samples of each
level:
\begin{equation}\label{eq:optimal_nsamples}
 N_{\ell} \propto
\sqrt{\frac{\V[Y_{\ell}]}{C_{\ell}}}\quad \ell=0,\dots,L,
\end{equation}
where $C_{\ell}$ is the cost of
computing one sample at level $\ell$.  We refer to \cite{cliffe2011multilevel,
giles2015multilevel}
for additional details.

The key idea that leads to computational savings is that fewer samples are necessary
to estimate $\E[Y_{\ell}]$ on finer levels, because 
$\V[Y_{\ell}]\rightarrow 0$ as
$h_{\ell} \rightarrow 0$ as long as $Q_h$ converges to $Q$ in expectation. 
The number of samples needed on coarser levels is still large, however samples are
less expensive to compute. 
This balancing of errors across the levels of the hierarchy leads to significant
improvements in computational time, while maintaining a desired level of
accuracy of the estimate.

\section{Forward model problem} \label{sec:darcy}
We are interested in the simulation of steady state groundwater
flow, governed by Darcy's law, in a porous medium where the 
permeability is not fully known.
We consider the mixed Darcy equations 
given by
\begin{equation}
  \begin{array}{lcr}
    \frac{1}{k(\bx,\omega)}\bq(\bx,\omega) + \nabla p(\bx,\omega) = 0 & \mbox{ in } D, \\
    \nabla \cdot \bq(\bx,\omega) = 0 & \mbox{ in } D,
  \end{array}\label{eq:darcy}
\end{equation}
with homogeneous Neumann boundary conditions $\bq \cdot \bn = 0$ on
$\Gamma_N$ and Dirichlet boundary conditions $p = p_D$ on
$\Gamma_D$. Here, $\Gamma_N \in \partial D$ and $\Gamma_D \in \partial
D$ are non overlapping partitions of $\partial D$, and $\bn$
denotes the unit normal vector to $\partial D$.
The uncertain permeability field $k(\bx,\omega)$ is modeled as a log-normal
random field such that $\log [k(\bx,\omega)]$ has a covariance function belonging
to the \Mat family, so that the pressure $p$ and Darcy flux $\bq$ are random fields as well.

We consider the discretization of \eqref{eq:darcy} with a log-normal permeability field
using the mixed finite element method~\cite{boffi2013mixed,fortin1991mixed}; 
this particular problem formulation has been analyzed
in~\cite{graham2016mixed}.
Assuming we are given an unstructured mesh $\mathcal{T}_h$ exactly 
covering $D$,
we consider solutions of the Darcy flux $\bq_h$ in the the lowest order Raviart--Thomas
finite element space denoted by $\bR_h \subset \bR:=H(\div,D)$,
and the pressure $p_h$ in the finite element space of piecewise constant function
denoted by $\Theta_{h} \subset \Theta := L^2(D)$.
Given an input realization
$k_h(\bx,\omega)$ as discussed in Section~\ref{sec:gaussian-rf}, 
the resulting discretized saddle-point problem  
can be written as 
\begin{equation}\label{eq:mixedDarcyDiscrete}
    \Lambda(k)_h \mathbf{X}_h :=
    \begin{bmatrix}
                M(k)_h & B_{h}^T \\
                B_{h} & 0
    \end{bmatrix}
    \begin{bmatrix}
        \bq_{h} \\ 
        \bp_{h}
    \end{bmatrix}
    =
    \begin{bmatrix}
        \f_{h} \\ 
        0
    \end{bmatrix}
    := \mathbf{G}_h,
\end{equation}
where $\f_h$ stems from the discretization of the Dirichlet boundary condition $p = p_D$ on $\Gamma_D$.

\subsection{Multilevel formulation}\label{subsec:ml-darcy}
We now consider the discretization of \eqref{eq:darcy} on a hierarchy of levels,
as needed by the MLMC algorithm.
Given the unstructured mesh $\T_h$ of $D$, we assume that 
a sequence of unstructured meshes  
$\T_{\ell}$ for $\ell=1,\dots,L$ on $D$ 
has been generated by recursively agglomerating
finer level elements. We denote the finest level mesh $\T_h$ as $\T_0$,
whereas $\T_L$ corresponds to the coarsest level.
For each coarse level, we construct the corresponding finite element spaces
$\bR_{\ell}$, $\Theta_{\ell}$, associated with the (agglomerated) mesh
$\T_{\ell}$ using methodology from AMGe methods, so that we are able to construct
operator-dependent coarse spaces for $H(\div)$ problems with guaranteed
approximation properties on general, unstructured grids;
see~\cite{lashuk2014construction,lashuk2012element,pasciak2008exact,
kalchev2016upscaling} for
further details.

We denote the piecewise constant interpolation operators from
coarser space $\Theta_{\ell+1}$ to the finer space $\Theta_{\ell}$ as
$P_{\theta}$ for $\ell = 0,\dots,L-1$.  
We also define the operators from the coarser space $\bR_{\ell+1}$ to
the finer space $\bR_{\ell}$ as $P_{\bu}$.
These operators are constructed
using techniques from AMGe; see~\cite{lashuk2014construction, lashuk2012element,
pasciak2008exact} for details about the operators
$P_{\theta}$ and $P_{\bu}$.

The discrete saddle-point block matrices are labeled ${\Lambda(k)}_{\ell}$, 
corresponding to the pair of finite element spaces $\bR_{\ell}$, $\Theta_{\ell}$
for $\ell=0,\dots,L$. Then, the discrete saddle-point problem \eqref{eq:mixedDarcyDiscrete}
at coarse level $\ell=1,\ldots,L$ reads

\begin{equation}\label{eq:ml-mixedDarcyDiscrete}
    \Lambda(k)_{\ell} \mathbf{X}_{\ell} :=
    \begin{bmatrix}
                M(k)_{\ell} & B_{\ell}^T \\
                B_{\ell} & 0
    \end{bmatrix}
    \begin{bmatrix}
        \bq_{\ell} \\
        \bp_{\ell}
    \end{bmatrix}
    =
    \begin{bmatrix}
        \f_{\ell} \\
        0
    \end{bmatrix}
    := \mathbf{G}_{\ell},
\end{equation}
where
\begin{equation}\label{eq:galerkin_projection_model}
M(k)_{\ell} := P_\bu^T M(k)_{\ell-1} P_\bu, \quad B_\ell:= P_\theta^T B_{\ell-1} P_\bu, \quad \mathbf{f}_{\ell} := P_{\bu}^T\mathbf{f}_{\ell-1}
\end{equation}
are the Galerkin projection at level $l=1,\dots,L$ of the corresponding fine grid matrices and vectors.

For an efficient MLMC simulation, it is necessary to repeatedly assemble and solve \eqref{eq:ml-mixedDarcyDiscrete}
on coarse levels $\ell=1,\dots,L$ for various realizations of $k_{\ell}(\bx,\omega^{(i)})$ \emph{without visiting the fine grid}.

Since $B_\ell$ and $\mathbf{G}_{\ell}$ are independent of $k_{\ell}(\bx,\omega^{(i)})$, such matrices and vectors can be computed once
--- using the Galerkin projection in \eqref{eq:galerkin_projection_model} --- before the MLMC simulation.
However, the efficient computation of $M(k)_{\ell}$ on coarse levels requires particular care, since this matrix depends on the random
parameter $k_{\ell}(\bx,\omega^{(i)})$.
To this aim, we exploit the sophisticated data structures of the AMGe hierarchies,
which closely mimic the same data structures of geometric multigrid and include
topological tables (i.e., element-element, element-face connectivity) and degree-of-freedom 
to element mapping for all levels of the hierarchy.
Specifically, at each level $\ell$  we denote with
$L^{\bR_{\ell}}_{e_{\ell}}$ the mapping between local (to the element $e_{\ell}\in \T_{\ell}$)
and global (for the space $\bR_{\ell}$) degrees of freedom.
These local to global mappings are then used to assemble local (to each agglomerated element
$e_{\ell} \in \T_{\ell}$ ) mass matrices into the global one. 
Then, given the piecewise-constant on the elements of $\T_{\ell}$ input $k_{\ell}$, 
the global weighted mass matrix $M(k)_{\ell}$ for the space $\bR_{\ell}$ is computed as 
\begin{equation}\label{M(k) element assembly}
M(k)_{\ell}=\sum_{e_{\ell} \in \T_{\ell}}s_{e_{\ell}}(L_{e_{\ell}}^{\bR_{\ell}})^T M_{e_{\ell}} (L_{e_{\ell}}^{\bR_{\ell}}) , 
\end{equation}
where $s_{e_{\ell}}=\left.k_\ell^{-1}\right|_{e_{\ell}}$ is the restriction of $k_\ell^{-1}$
to the element $e_{\ell}\in\T_{\ell}$,
and $\curly{M_{e_{\ell}}}_{e_{\ell}\in \T_{\ell}}$ are the 
local mass matrices for the space $\bR_{\ell}$.
These local matrices are computed once during the construction of the AMGe hierarchy by local 
(to each agglomerated element) Galerkin projection of partially assembled mass matrices from 
the previous (finer) level, and then stored for future use in the MLMC simulation.
For details on the assembly procedure see \cite{ChristensenVillaVassileski15,ChristensenVillaEngsigEtAl17},
where a time-dependent two-phase porous media flow is solved with optimal complexity
on coarse (upscaled) levels, and \cite{ChristensenVassilevskiVilla17}, where a nonlinear scalable multilevel solver
for single-phase porous media flow is presented.

\subsection{Linear solution of forward model problem}
\label{subsec:darcy-solver}
The linear system \eqref{eq:ml-mixedDarcyDiscrete} on each level 
is iteratively solved using preconditioned GMRES.
We consider a preconditioner based on the 
approximate block-LDU factorization of the operator ${\Lambda(k)}_{\ell}$, given by:
\begin{equation}
  \label{eq:ldu}
  {\M(k)}_{\ell} =
  \begin{bmatrix}
    I & \\
    {B}_{\ell}\tilde{M}(k)_{\ell}^{-1} & I
  \end{bmatrix}
  \begin{bmatrix}
    \tilde{M}(k)_{\ell} & \\
    & -\tilde{S}_{\ell}
  \end{bmatrix}
  \begin{bmatrix}
    I & \tilde{M}(k)_{\ell}^{-1}{B}_{\ell}^{T}\\
    & I
  \end{bmatrix}.
\end{equation}
Here, $\tilde{M}(k)_{\ell}$ is a \emph{cheap} preconditioner for the
mass matrix, ${M(k)}_{\ell}$, such as a diagonal approximation or a
small number of Gauss-Seidel iterations. We use three Gauss-Seidel
iterations in our numerical experiments.
$\tilde{S}_{\ell} = B_{\ell}\diag(M(k)_{\ell})^{-1}B_{\ell}^T$ is the
approximate Schur-complement, which is symmetric positive definite and
sparse.  In our numerical experiments, we approximate the action of
the approximate Schur-complement inverse by a single algebraic multigrid
V-cycle; specifically, we use BoomerAMG from the solvers library
\textit{hypre}~\cite{hypre}.

\begin{remark}
  It is worth noticing that the dominant cost in applying the preconditioner \eqref{eq:ldu}
  is to approximate the action of $\tilde{S}_{\ell}^{-1}$.
  This justifies the use of a full LDU factorization instead of simpler 
  methods, such as block-diagonal or even block-triangular approaches~\cite{MFMurphy_GHGolub_AJWathen_2000a}.
  Our numerical studies showed that the performance of the full LDU approach was better not only in terms of
  number of iterations, but also in total solve time compared to the other methods.
  \label{remark:LDU}
\end{remark}

\section{Gaussian Markov random field based sampling techniques for spatially correlated random fields} \label{sec:gaussian-rf}
In this section we discuss generating
realizations of a log-normal random
field,
$k(\bx,\omega) = \exp[\theta(\bx,\omega)]$, where $\theta(\bx,\omega)$ is a
Gaussian random field with a certain mean and covariance structure to be
used as input coefficients for a MLMC simulation. 
We extend the sampling strategy of \cite{o16} to include a scalable domain embedding technique
allowing the use of non-matching meshes to sample from the Gaussian random field $\theta(\bx,\omega)$. 

In particular, we consider the \Mat family of covariance functions, which is a
common choice in geostatistics~\cite{chiles2009geostatistics}.
The \Mat covariance function is given by 
\begin{equation}
    \operatorname{cov}(\bx, \by) = \frac{\sigma^{2}}{2^{\nu -1}\Gamma(\nu)} (\kappa \Vert
    \by-\bx\Vert)^{\nu}K_{\nu}(\kappa \Vert \by-\bx\Vert),
    \label{eq:matern}
\end{equation}
where $\sigma^{2}$ is the marginal variance, $\nu>0$ determines the mean-square
differentiability of the underlying process, $\kappa>0$ is a scaling factor
inversely proportional to the correlation length, $\Gamma(\nu)$ is the gamma function,
and $K_{\nu}$ is the modified
Bessel function of the second kind~\cite{matern1986spatial}.

To realize a sample of a Gaussian random field with \Mat covariance as in \cite{o16}, 
we employ a sampling method 
that uses a link between Gaussian fields and Gaussian
Markov random fields established in \cite{lindgren2011explicit}.
The method utilizes the fact that the solution,
$\theta(\bx,\omega)$, of the fractional SPDE
given by 
\begin{equation} 
    (\kappa^{2} - \Delta)^{\alpha/2}\theta(\bx,\omega) = g
    \W(\bx,\omega) \quad  \bx \in \R^d(d=2\mbox{ or }3), \ \alpha = \nu +\frac{d}{2}, \
\kappa> 0, \nu >0, \label{eq:spde_fractional} 
\end{equation} is a Gaussian field with
underlying \Mat covariance~\cite{whittle1954stationary,whittle1963stochastic}. Above, $\W$ is Gaussian white noise, and the scaling
factor $g$ is chosen to impose unit
marginal variance of the random field as 
\begin{equation*} 
    g = (4\pi)^{d/4} \kappa^{\nu}  \sqrt{\frac{\Gamma\left(\nu +
    d/2\right)}{\Gamma(\nu)}}.
\end{equation*} 
We restrict the smoothness parameter, $\nu$, to be of the form $\nu=\alpha - d/2$ for an
even integer $\alpha$. 
Specifically, in three dimensions, 
the choice of $\nu=1/2$ results in the random field $\theta(\bx,\omega)$
having an underlying exponential covariance structure, as \eqref{eq:matern}
reduces to $\operatorname{cov}(\bx,\by) = \sigma^2 e ^{-\kappa \norm{\by - \bx}}.$

Then, \eqref{eq:spde_fractional} 
reduces to the following standard reaction-diffusion equation: 
\begin{equation} 
    (\kappa^{2}-\Delta)\theta(\bx,
    \omega) = g\mathcal{W}(\bx,\omega) \quad \bx\in\R^d, \kappa>0.
\label{eq:SPDE} 
\end{equation} 
Thus a scalable sampling method is
equivalent to efficiently solving the stochastic reaction-diffusion equation
given by \eqref{eq:SPDE}. 

It should be noted that defining covariance operators as fractional inverse
powers of differential operators is a common approach for the solution of
large-scale Bayesian inverse problems governed by PDE forward models, see e.g.,
\cite{Stuart10,Bui-ThanhGhattasMartinEtAl13}, as it allows for efficient
evaluation of the covariance operator using a fast and scalable multigrid
solver.  For further details on the approximation of Gaussian random fields with
\Mat covariance functions using the Gaussian Markov random field representation
of a SPDE, we refer to~\cite{lindgren2011explicit,
simpson2012order, lindgren2015bayesian}.

\subsection{Stochastic PDE sampler}
We consider the solution of the stochastic reaction-diffusion
equation given by \eqref{eq:SPDE} to produce realizations of a Gaussian random
field using the mixed finite element method on a bounded domain $D \subset \R^d$.
When posing the SPDE on a bounded domain, boundary conditions must be imposed,
however the proper boundary conditions for the stochastic fields are an open research
problem; see, e.g., \cite{DaonStadler2016}. We consider using deterministic homogeneous Neumann boundary conditions (zero
normal-derivatives); however, this choice introduces boundary artifacts that inflate the variance along the boundary of the
domain, as observed in~\cite{lindgren2011explicit}. 

One approach to mitigate this issue is to extend the domain
of interest by a distance greater than the correlation length, 
solve \eqref{eq:SPDE} on the
extended domain, then restrict the solution back to the original domain to
generate a Gaussian field realization. This procedure was explored in
\cite{o16}, and mitigates the artificially inflated variance 
as the boundary effects are negligible at
a distance greater than the correlation length away from the 
boundary~\cite{lindgren2015bayesian}.
However, this can pose a challenge 
with complicated domains and/or unstructured meshes.   
In this work, we propose an alternative domain embedding technique with two 
non-matching meshes, using
a scalable transfer of discrete fields between the two meshes which can be
arbitrarily distributed among different processors, see Section \ref{sec:projection}. 
First the SPDE is
discretized using the mixed finite element method on an extended regular domain,
then the finite element solution is transferred to the original domain resulting
in a realization of the Gaussian random field for use as the input coefficient
in a Monte Carlo simulation.

\subsection{Mixed finite element discretization}
Let $D$ be a given polygonal/polyhedral domain with an unstructured
mesh $\T_h$. We embed $D$ in an extended regular domain (e.g., a box) ${\oD}$
meshed by ${\oT}_{\oh}$, where ${\oh}
\simeq h$. In contrast to the approach taken in \cite{o16}, the two meshes here do not necessarily match on $D$; in addition,
${\oT}_{\oh}$ does not necessarily respect $\partial D$. We
assume that ${\oT}_{\oh}$ is obtained by several steps of
refinement of an initial coarse mesh.
We consider the mixed finite element
discretization~\cite{boffi2013mixed,fortin1991mixed} of \eqref{eq:SPDE}
on the regular domain ${\oD}$ with mesh ${\overline
\T}_{\overline h}$.
We introduce the functional spaces ${\oR}:= H(\div,\oD)$ and 
${\oTheta}:= L^2(\oD)$, as well as the bilinear forms 
\begin{equation*}
    \begin{array}{lll}
        m(\obu, \obv)  &:=(\obu, \obv)       & \forall \, \obu, \obv \in \oR,\\
        w(\otheta, \oq) &:= (\otheta, \oq)     & \forall \, \otheta, \oq \in \oTheta,\\
        b(\obu, \oq)    &:= (\div\, \obu, \oq) & \forall \, \obu \in \oR, \oq \in
        \oTheta,\\
    \end{array}
\end{equation*}
and the linear form
\begin{equation*}
    \oF^\W(\oq) := (\W, \oq) \quad \forall \, \oq \in \oTheta.
\end{equation*}
Above the symbol $(\cdot, \cdot)$ denotes the usual inner product
between scalar (vectorial) functions in $L^2(D)$ ($[L^2(D)]^d$). 

Let $\oR_{\oh} \subset \oR$ denote the lowest order Raviart--Thomas space and
$\oTheta_{\oh} \subset \oTheta$ denote the finite element space of piecewise
constant functions defined
on the fine triangulation $\oT_{\oh}$ of $\oD$, then   
we seek the solution of the mixed finite
element discretization given by 
\begin{problem}
  Find $\obu_{\oh} \in {\oR}_{\oh}$ and
  ${\otheta}_{\oh} \in {\oTheta}_{\oh}$ such that 
  \begin{equation}
    \begin{array}{ll} 
        m(\obu_{\oh},\;\obv_{\oh})+ b(\obv_{\oh},\;\otheta_{\oh})= 0 & \mbox{ for
        all } \obv_{\oh} \in \oR_{\oh},\\ 
      b(\obu_{\oh},\;\oq_{\oh})-\kappa^2\;w(\otheta_{\oh},\;\oq_{\oh})  =
      -g\oF^\W(\oq_{\oh})& \mbox{ for all } \oq_{\oh} \in  {\oTheta}_{\oh}  
    \end{array} \label{eq:reac_diff_weak}
  \end{equation}
  with essential boundary conditions $\obu_{\oh} \cdot \bn = 0$ on $\partial {\overline D}$.  
\end{problem}

To formulate the linear algebra representation of the stochastic right
hand side $\oF^\W(\oq_{\oh})$, two properties of Gaussian white noise defined on
a domain $\oD$ are used. 

For any set of test functions $\curly{\oq_i \in L^2(\oD), i=1,\dots,n},$
the expectation and covariance measures are given by
\begin{align}
\mathbb{E}[ (\oq_{i},\mathcal{W})] &= 0, \\
\operatorname{cov}\left( (\oq_{i},\mathcal{W}), (\oq_{j}, \mathcal{W})\right)
&= (\oq_{i},\oq_{j}).
\end{align}
By taking $\oq_{i},\oq_{j}$ as piecewise constants so that $\oq_{i},\oq_{j}\in
\oTheta_{\oh}$, the
second equation implies that the covariance measure over a region of the domain
is equal to the volume of that region \cite{lindgren2011explicit}.
Then the discrete stochastic linear right hand side is given by
\begin{equation*}
    {\overline f}_{\oh} = {\overline W}_{\oh}^{\frac{1}{2}}\xi_{\oh}(\omega),
    \quad \xi_{\oh}(\omega) \sim \mathcal{N}(0,I)
\end{equation*}
where ${\overline W}_{\oh}$ is the mass matrix for the space
$\oTheta_{\oh}$, and $\mathcal{N}(0,I)$ denotes the multivariate normal distribution
with zero mean and covariance matrix $I$, that is, each component of
$\xi_{\oh}(\omega)$ is standard normal and the components are independent.
It should be noted that the mass matrix ${\overline W}_{\oh}$ for the space
$\oTheta_{\oh}$ is diagonal, hence its square root can be computed cheaply.

Then the discrete mixed finite element problem can be written as the linear system
\begin{equation}
{\overline \A}_{\oh}{\overline U}_{\oh} =
{\overline F}_{\oh},
\label{eq:linsys}
\end{equation}
with block matrix and block vectors defined as
\begin{equation}
  {\overline \A}_{\oh} =
  \begin{bmatrix}
    {\overline M}_{\oh} & {\overline B}_{\oh}^T \\
    {\overline B}_{\oh} & -\kappa^2 {\overline W}_{\oh}
  \end{bmatrix},
  \quad
  {\overline U}_{\oh} =
  \begin{bmatrix} \obu_{\oh}  \\ \obtheta_{\oh} \end{bmatrix},
  \quad
  {\overline F}_{\oh} =
  \begin{bmatrix} 
    0  \\ -g\,{\overline f}_{\oh}(\omega)
  \end{bmatrix},
  \label{eq:linsys_blocks}
\end{equation}
where ${\overline f}_{\oh}(\omega) \sim \mathcal{N}(0,{\overline W}_{\oh}),$
${\overline M}_{\oh}$ is the mass matrix for the space $\oR_{\oh}$,
${\overline B}_{\oh}$ stems from the discretization of the divergence operator,
${\overline W}_{\oh}$ is the (diagonal) mass matrix for the space
$\oTheta_{\oh}$, and $\obu_\oh$, $\obtheta_{\oh}$ are the coefficient vectors
of the finite element functions
when expanded in terms of the respective basis functions.
We remark that the covariance structure of the samples generated by solving the mixed form \eqref{eq:SPDE}
is equivalent to the covariance structure of samples obtained by solving the primal form of the SPDE in \cite{lindgren2011explicit},
as shown in \cite{o16}. In fact, using simple algebraic manipulation, it is immediate to show that $\obtheta_{\oh} \sim \mathcal{N}( 0, {\overline{C}}_{\oh} )$, where $\overline{C}_{\oh} = {\overline S}_{\oh}^{-1}{\overline W}_{\oh}{\overline S}_{\oh}^{-T}$ and ${\overline S}_{\oh} = \kappa^2 {\overline W}_{\oh} + {\overline B}_{\oh}{\overline M}_{\oh}^{-1}{\overline B}_{\oh}^T$ stems from a nonlocal
{\em discontinuous} Galerkin (interior penalty) discretization of the original PDE \eqref{eq:SPDE}; cf. \cite{RVW}.

\section{$L^2$-Projection} \label{sec:projection}
This section concerns a parallel scalable $L^2$-projection of discrete fields between non-matching
partially overlapping meshes. Specifically, in the MLMC simulation we use the $L^2$-projection to
transfer a realization of a Gaussian random field $\otheta_\oh$ on ${\overline
\T}_{\oh}$ --- computed by solving \eqref{eq:linsys} approximately by an iterative method ---
to the unstructured mesh $\T_{h}$ of the original domain $D$, where we then solve the forward problem
\eqref{eq:mixedDarcyDiscrete}.

\subsection{Projection between non-matching meshes} \label{sec:l2_projection}

Let $\Theta_h = \text{span}\;\{\varphi_\tau\}_{\tau \in \T_h}$ and 
$\oTheta_\oh = \text{span}\;\{{\ovarphi}_{\otau}\}_{\otau \in \oT_\oh}$, 
and assume we have the functions $s \in \Theta_h$ and $\os \in
\oTheta_\oh$. Writing the functions in terms of their respective bases, we have 
$s=\sum_{\tau} s_{\tau}\varphi_\tau$ and
$\os=\sum_{\otau} \os_{\otau}{\ovarphi}_{\otau}$, and computing the quantity 
$(s,\;{\ovarphi}_{\otau})$ yields
\begin{equation*} 
    (s,\;{\overline \varphi}_{\overline \tau}) =
    \sum\limits_{\tau \in \T_h} s_\tau 
    \int\limits_{\tau \cap {\overline \tau}}
    \varphi_\tau {\overline \varphi}_{\overline \tau}\;d\bx \mbox{  for
    all } {\overline \tau} \in {\overline \T}_\oh.  
\end{equation*} 

Introducing the matrix ${\overline G} = (g_{{\otau},\tau})$ where
$g_{{\otau},\tau} = \int\limits_{\tau \cap {\otau}}\varphi_\tau {\ovarphi}_{\otau}\;d\bx,$ 
we can rewrite the integral moments in matrix-vector
form as follows
\begin{equation}\label{matrix vector product with overline G}    
    {\overline G} \bm{s} = \left (\sum\limits_{\tau} g_{{\overline \tau},\tau} s_\tau \right
      )_{{\overline
      \tau} \in {\overline \T}_\oh},\;\; \bm{s} = ( s_\tau )_{\tau \in \T_h}.
\end{equation}

We are interested in the $L^2$-projection $s = \sum_{\tau}
s_{\tau}\varphi_\tau$ of $\os$ onto $\Theta_h$.
Since $\os=\sum_{\overline \tau} \os_{\overline \tau}{\overline
\varphi}_{\overline \tau}$, 
we have 
\begin{equation*} 
    (\os,\;{\varphi}_{\tau}) =
    \sum\limits_{\overline \tau \in \oT_\oh} \os_\tau 
    \int\limits_{\tau \cap {\overline \tau}}
    \varphi_\tau {\overline \varphi}_{\overline \tau}\;d\bx \mbox{  for
    all } {\overline \tau} \in {\overline \T}_\oh,   
\end{equation*} 
which can be written in matrix-vector
form as $W_h \bm{s} = {\overline G}^T \bm{\os}$ where $\bm{\os}$ is the vector of
coefficients with entries $\os_{\otau}$. Therefore, letting $G = {\overline G}^T$ we have 
\begin{equation}
   \bm{s} = W_h^{-1}{\overline G}^T\bm{\os} = W_h^{-1} G \bm{\os} := \Pi_h \bm{\os},  
\label{eq:proj_op}
\end{equation}
where $W_h$ is the (diagonal) mass matrix for the
space $\Theta_h$.

\subsection{Implementation details of the $L^2$-projection operator}
\label{sub-sec:mapping}
The Petrov-Galerkin assembly of the discrete $L^2$-projection $\Pi_h$ 
requires computing the intersection between elements, and building a suitable set of quadrature points and weights from the intersection.
Given two meshes $\T_h$ and ${\overline \T}_h$
we search for each pair of elements $\tau \in \T_h$ and
${\overline \tau} \in {\overline \T}_{\overline h}$ with intersection
$I = \tau \cap {\overline \tau} \neq \emptyset$.
We mesh the intersection $I$ into a simplicial complex $I_h$ and map a suitable quadrature rule, such as standard Gaussian formulas~\cite{stroud1966gaussian}, to each simplex $S \in I_h$.
We transform the resulting quadrature points to the reference configuration of
both elements ${\tau}$ and ${\overline \tau}$ and perform the assembly procedure.

One critical aspect for run time performance is intersection detection~\cite{ericson2004rcd}.
Of particular interest are linear time complexity algorithms such as spatial hashing~\cite{lefebvre2006psh} for quasi-uniform meshes, or advancing front algorithms~\cite{gander2009ngp} for meshes with varying size elements.

In large-scale parallel computations, meshes are generally arbitrarily
distributed and no relationship between their elements is explicitly available.
For determining such relationships based on spatial information, we consider parallel
algorithms relying on space filling curves~\cite{bader2012space} and parallel
tree searches~\cite{krause2016parallel}. We have implemented the algorithm
in~\cite{krause2016parallel} by exploiting the software libraries
\MFEM{}~\cite{mfem}, for handling the finite element representations, and
\MOONOLITH{}~\cite{moonolith}, 
for handling the parallel intersection detection/computation and automatic load-balancing.

Here we summarize the overall parallel search approach which we implemented with MPI~\cite{mpi-3.0}.
Our strategy exploits the implicit and
self-affine structure of octrees for adaptively constructing a bounding volume
hierarchy (BVH) fitting the volume of interest. 
While the octree is a global object, hence the volume described by the root cell $C_0 \supseteq D \cup \oD$
is the same for all processes, and can be refined in the same way by any process without any communication,  
the BVH is constructed using local geometric information that needs to be exchanged for intersection testing.
Similar to the octree, the bounding volumes associated with the BVH nodes are axis-aligned bounding-boxes (AABBs).
However, instead of following the octree predefined subdivision pattern, the AABBs of the BVH fit the data more accurately and do not
necessarily form a partition.
For each octree cell $C_n$ of node $n$, we have an associated bounding-box $B_n^p$ which is part of the BVH constructed by process $p$.
The extra BVH allows us to have tighter bounding volumes for the nodes at coarse levels of our hybrid hierarchy, hence allowing early pruning when performing tree-searches and dramatically improving the performance of the search. 
The goal of using this hybrid octree/BVH data-structure is to perform a cheap broad-phase intersection test without the need
of exchanging mesh data.

Next we give a brief explanation of the several steps of our algorithm which are analyzed in terms of parallel performance in Section~\ref{sec:results}, Figures \ref{fig:crooked_pipe_sampler_setup} and \ref{fig:SPE10_sampler_setup}. These steps consists of \emph{element bounding volumes generation}, \emph{BVH comparison}, \emph{load balancing}, \emph{matching and rebalancing}, and \emph{computation} of the $L^2$-projection operator. 

In the \emph{element bounding volumes generation} step, we compute the AABB 
for each element of the input meshes which we use
for inserting the element in the octree/BVH data-structure.

In the \emph{BVH comparison} step, we construct both the octree/BVH and we perform a search in the branches of the tree where there are potential intersections.
In a potential intersection region containing $D \cap \oD$ we adaptively refine the tree in an iterative fashion.
At each iteration we refine the octree/BVH, and exchange the necessary information for constructing the search paths (or tree-traversals) consistently.
Each pair of processes $\{p, q\}$ has dedicated search-paths  which are updated/refined only when the required data (MPI-message) is available using asynchronous point-to-point communication. This data consists of the bounding-volumes $B_n^p, B_n^q$, which are tested for intersection, and the number of elements associated with each octree/BVH node $n$. If $B_n^p \cap B_n^q = \emptyset$ and if either $p$ or $q$ do not have elements associated with $n$, we stop the search for the sub-tree with root $n$ for the pair of processes $\{p, q\}$. 

From the BVH comparison we obtain a list $L$ of tuples $\{n, p, q\}$, where $n$ is a node of the octree/BVH and $p$ is a process
having elements intersecting with the bounding-box $B_n^q$ of process $q$. Note that if $\{n, p, q\}$ exists then $\{n, q, p\}$ also exists and we consider them to be the same tuple.
For any pair of processes $p, q$ having an entry in $L$ for node $n$ we estimate the cost of performing the intersection test
between the sets of elements associated with $n$ by a cost function $\gamma_n(s_p, s_q) = s_p  s_q$, where $s_k, k \in \{p, q\}$, is the size of the set of elements associated with $n$ in the memory of process $k$.

The latter step allows us to perform the \emph{load balancing} task. 
The load balancing is done by assigning the elements associated with each tuple $\{n, p, q\}$ in such a way that the work, according to the cost function $\gamma_n(s_p, s_q)$, is distributed as evenly as possible among processes. 
The load balancing algorithm exploits the ordering computed by linearization of the octree (i.e., Morton ordering) for splitting the work and grouping together nearby elements. 
Note that potential imbalances which might be caused by the output-sensitivity (i.e., the cost of the computation is influenced by the size of the output) of the problem due to both position and distribution of the elements are mitigated by the search procedure. 
This is only feasible because we delay the element-to-element intersection test to the latest possible moment. 
In fact, the next step which is the \emph{matching and rebalancing} step, consists of communicating the necessary elements and determining the
intersection pairs. Once we have the intersection pairs we re-balance once again, but this time at much finer granularity, to ensure an efficient \emph{computation} step consisting on the actual element-to-element intersection computation and numerical quadrature for the assembly of the $L^2$-projection operator.

Our intersection-detection approach enables an efficient broad-phase intersection testing when the two (or more) finite element meshes are partially overlapping, which is the case for the problem presented in this paper.
However, for the case where we have the prior knowledge that the two meshes are describing the same volume refining the tree before the construction
of the search paths might provide slight performance improvements in tree-search phase of the intersection detection algorithm.
For more details on parallel variational transfer we refer to~\cite{krause2016parallel}.

\section{Hierarchical SPDE Sampler} \label{sec:ml-sampler}
In this section we describe our proposed hierarchical SPDE sampling technique
using domain embedding with non-matching meshes for MLMC simulations. 
We first describe the process of
generating the sequences of coarser levels of $\oT_\oh$, 
introducing the necessary finite element spaces and
interlevel operators that will be used, followed by 
implementation details of mapping 
between non-matching meshes for coarse levels. 
Finally, the iterative solution strategy of the saddle-point problem 
is described for all levels.

\subsection{Multilevel structure}\label{subsec:ml} 
We now describe the multilevel structure of the SPDE sampler using 
the structured grid of the regular domain $\oD$. 
We assume that we have a nested sequence 
of meshes ${\oT}_{\ell}$ on the regular domain ${\oD}$ for $\ell=0,\dots,L$ obtained by several
steps of refinement of an initial coarse mesh $\oT_{h_L}$ with $\oT_0=\oT_{\oh}$. 
The corresponding finite element spaces $\oR_{\ell}$, $\oTheta_{\ell}$ are 
constructed for each coarse level to form a geometric hierarchy of
standard Raviart--Thomas finite element spaces.

For $\ell = 0, \dots, L$, we denote the saddle-point block matrices
$\overline{\A}_{\ell}$ corresponding to the pair of finite element spaces $\oR_{\ell}$,
$\oTheta_{\ell}$.
The block interpolation operator for the matrix 
$\overline{\A}_{\ell}$ is defined as
\begin{equation}\label{eq:mono_mg_interp}
  \overline{\P} =
  \begin{bmatrix} 
    \overline{P}_{\bu} & 0 \\
    0 & \overline{P}_{\theta} 
  \end{bmatrix},
\end{equation}
where 
$P_{\bu}$ is the operator from the coarser space $\oR_{\ell+1}$ to
the finer space $\oR_{\ell}$, and
$\oP_{\theta}$ is the piecewise constant operator
from coarser space $\oTheta_{\ell+1}$ to the finer space $\oTheta_{\ell}$
for $\ell=0,\dots,L-1$.
These interpolation operators for the structured hierarchy of
uniformly refined meshes
are the canonical interpolation operators of geometric
multigrid.
Then, we write \eqref{eq:linsys}  
at coarse level $\ell=1,\dots,L$ as
\begin{equation}
{\overline \A}_{\ell}{\overline U}_{\ell} =
{\overline F}_{\ell},
\label{eq:ml-linsys}
\end{equation}
where 
\begin{equation*}
\overline{\A}_{\ell} := \overline{\P}^T \overline{\A}_{\ell-1}
\overline{\P}, \quad  {\overline F}_{\ell} := \overline{\P}^T {\overline F}_{\ell-1}.
\end{equation*}
To conclude this section, we present the multilevel definition of the
$L^2$-projection operator ${\overline G}_0 = {\overline G}$ presented in Section~\ref{sec:l2_projection}, see formula \eqref{matrix vector product with overline G}.       
By letting $\Pi_0$ denote $L^2$-projection operator on the fine mesh and $G_0 = {\overline G}^T_0$ denote the Petrov-Galerkin mass operator between the non-matching mesh at the fine grid level, we recursively define the $L^2$-projection operator between coarse meshes as
\begin{equation}\label{eq:recursive_l2_projection}
\Pi_{\ell+1} = W_{\ell+1}^{-1} G_{\ell+1}, \quad G_{\ell+1} = P_{\theta}^T G_{\ell}\oP_{\theta},
\end{equation}
where $W_{\ell+1}$ denotes the mass matrix in the space $\Theta_{\ell+1}$,
and $P_{\theta}$ is the unstructured hierarchy's interpolation operator 
discussed in Section~\ref{subsec:ml-darcy}.

\subsection{Hierarchical SPDE sampler with non-matching mesh embedding} 
We have shown that a realization of a Gaussian random field $\theta_h$ on $\T_h$ can be obtained by
solving the linear system \eqref{eq:linsys} for $\obtheta_{\oh}$, then computing
\begin{equation} 
    {\btheta}_h = \Pi_h {\obtheta}_\oh. 
    \label{eq:theta_l2projection}
\end{equation} 
From the linearity of the $L^2$-projection, it immediately follows that
${\btheta}_h \sim \mathcal{N}({\bf 0}, C_h)$, with $C_h = \Pi_h \overline{C}_{\oh}
\Pi_h^T$.
This method proves to be scalable and efficient as it is able to leverage existing
solvers and preconditioners for saddle-point problems with 
structured grids; however, the parametrization
of $\theta$ is mesh-dependent. 
For MLMC, a realization of a Gaussian
random field must be computed on a fine and coarse spatial resolution for the
same random event $\omega$. 
Thus, we consider generating $\theta_{\ell}(\omega)$ and $\theta_{\ell+1}(\omega)$
for the same $\omega$.

As shown in~\cite{o16}, the Gaussian random field $\otheta_{\ell}(\omega)$ whose
coefficient vector is given by
\begin{equation}\label{eq:ml_spde_sampler}
 \begin{bmatrix} 
      {\overline \bu}_{\ell} \\ \obtheta_{\ell}(\omega) 
    \end{bmatrix} = \overline{\mathcal{A}}^{-1}_{\ell} 
    \begin{bmatrix} 0 \\ -g{\overline W}^{1/2}_{\ell} \xi_{\ell}(\omega)\ \end{bmatrix}
\end{equation}
admits the following two-level decomposition:
\begin{equation}\label{eq:bileveld}
    \obtheta_{\ell}(\omega)={\overline
    P}_{\theta}\obtheta_{\ell+1}(\omega)+\delta\obtheta_{\ell}(\omega),
\end{equation}
where $\obtheta_{\ell+1}(\omega)$ is a coarse representation of a Gaussian random field from
the same distribution on $\oTheta_{\ell+1}$, and
\begin{equation}
  \begin{bmatrix} 
    \overline{\mathcal{A}}_{\ell} & \overline{\mathcal{A}}_{\ell} \overline{\mathcal{P}} \\ 
    \overline{\mathcal{P}}^T \overline{\mathcal{A}}_{\ell} & 0 
  \end{bmatrix} 
  \begin{bmatrix} 
    \delta \overline{U}_{\ell}  \\ 
    \overline{U}_{\ell+1}
  \end{bmatrix} 
= \begin{bmatrix} 
    \overline{\mathcal{F}}_{\ell}  \\ 0
  \end{bmatrix},
  \label{eq:ml_sampling}
\end{equation}
with the block expressions given by
$$ \delta \overline{U}_{\ell} = 
  \begin{bmatrix} \delta \overline{\bu}_{\ell} \\ \delta \obtheta_{\ell}(\omega) \end{bmatrix}, 
  \; \overline{U}_{\ell+1}=  
  \begin{bmatrix} \overline{\bu}_{\ell+1} \\  \obtheta_{\ell+1}(\omega) \end{bmatrix}, 
  \mbox{ and } \overline{\mathcal{F}}_{\ell}= 
  \begin{bmatrix} 0 \\ -g \overline{W}_{\ell}^{1/2}\xi_{\ell}(\omega) \end{bmatrix}. $$

Given $\xi_{\ell}(\omega)$, we compute the realizations of the Gaussian field 
at levels $\ell$ and $\ell+1$ on the spaces $\Theta_{\ell}$ and $\Theta_{\ell+1}$ 
by first computing $\obtheta_{\ell+1}$ by solving the saddle-point system
\eqref{eq:ml-linsys} at level $\ell+1$
using the methodology described in Section~\ref{subsec:solvers}.
Then we compute $\obtheta_{\ell}$ by iteratively solving
\eqref{eq:ml-linsys} at level $\ell$ with $\overline{\mathcal{P}} \overline{U}_{\ell+1}$ as the initial guess.
Finally, using the $L^2$-projection operators recursively defined by \eqref{eq:recursive_l2_projection}, 
we simultaneously transfer both $\otheta_{\ell}$ and $\otheta_{\ell+1}$ to the
unstructured mesh hierarchy,
and we write
\begin{equation*}
\btheta_{\ell} = \Pi_{\ell} \obtheta_{\ell}, \quad \btheta_{\ell+1} = \Pi_{\ell+1} \obtheta_{\ell+1}.
\end{equation*}

\subsection{SPDE sampler saddle-point problem linear solution}
\label{subsec:solvers}
We now discuss our methodology for the solution of the saddle-point system \eqref{eq:linsys}
using a scalable solver for $H(\div)$ problems.
The hybridization solver that we employ reduces the original saddle-point system to a 
symmetric, positive definite system, which after appropriate diagonal rescaling 
(cf. \cite{lee2017parallel}) is successfully solved by classical AMG solvers designed 
for $H^1$ equivalent problems.
The hybridization approach is a classical technique used for solving saddle-point 
problems arising from discretizations of mixed systems posed in $H(\div)$.
More specifically, in our setting, we have a saddle-point matrix of the form
\begin{equation*}
\left [
\begin{array}{cc}
{\overline M} & {\overline B}^T\\
{\overline B} & -{\overline W}
\end{array} \right ].
\end{equation*}

{\em Hybridization} refers to decoupling the degrees of freedom associated with the interfaces between the elements of the mesh corresponding to 
the first, vector unknown $\obu$ (coming from the Raviart--Thomas space), and then
imposing the difference of the decoupled quantities from both sides of the element interfaces to be zero posed as constraints using Lagrange multipliers.
In this way, one ends up with an equivalent saddle-point system with one extra set of unknowns, namely the Lagrange multipliers $\blambda$. The embedding system consists of the now decoupled (element-by-element) vector unknown
${\widehat \bu}$, the original piecewise constant unknowns $\obtheta$ and the Lagrange multipliers $\blambda$. The resulting matrix is symmetric
with saddle-point form
\begin{equation*}
\left [
\begin{array}{ccc}
{\widehat M} & {\widehat B}^T & {\overline C}^T\\
{\widehat B} & -{\overline W} & 0\\
{\overline C}            & 0 & 0
\end{array} \right ].
\end{equation*}
Above, ${\overline C}$ is the matrix coming from the constraint of zero jumps of the
decoupled unknowns ${\widehat \bu}$ across element interfaces.
The main property of the embedding matrix is that its two-by-two principal submatrix
\begin{equation*}
\left [
\begin{array}{cc}
{\widehat M} & {\widehat B}^T\\
{\widehat B} & -{\overline W}
\end{array} \right ]
\end{equation*}
is block-diagonal with blocks corresponding to degrees of freedom within each element (which are decoupled from the other elements).
Therefore, the reduced Schur-complement matrix
\begin{equation*}
\left [{\overline C},\;  0\right ] \left [
\begin{array}{cc}
{\widehat M} & {\widehat B}^T\\
{\widehat B} & - {\overline W}
\end{array} \right ]^{-1} \left [{\overline C},\;  0\right ]^T = 
{\overline C} \left ({\widehat M} + {\widehat B}^T {\overline W}^{-1}{\widehat B} \right )^{-1} {\overline C}^T,
\end{equation*}
is s.p.d.; it is explicitly available and provably equivalent to an $H^1$-discretization matrix. One problem (discussed and resolved in \cite{lee2017parallel}) is that depending on the choice of basis in the Raviart--Thomas space, one may need to diagonally rescale the Schur-complement so that the constant coefficient vector
corresponds to the constant function. The latter affects the successful use of
classical AMG methods (which implicitly assume that the constant vector is in the near null-space of the underlined matrix).

In our experiments the resulting hybridization linear system is solved with the conjugate gradient method preconditioned with \textit{hypre}'s highly scalable BoomerAMG solver~\cite{hypre}.

\section{Numerical Results} \label{sec:results}

We now demonstrate the numerical performance and parallel
scalability of the hierarchical SPDE
with non-matching mesh embedding, and include standard
results from MLMC computations using our proposed SPDE sampler 
for two different three-dimensional spatial geometries.

\begin{figure}[htbp]
\clearpage{}\tikzstyle{postprocess_block} = [rectangle, draw, fill=black!10,
    text width=5em, text badly centered, node distance=3cm, rounded corners, font=\footnotesize]
\tikzstyle{sampler_block} = [rectangle, draw, fill=red!20,
    text width=9em, text centered, node distance=1cm, rounded corners, minimum height=4em, font=\small]
\tikzstyle{projection_block} = [rectangle, draw, fill=blue!20,
    text width=8em, text centered, node distance=1cm, rounded corners, minimum height=4em, font=\small]
\tikzstyle{model_block} = [rectangle, draw, fill=green!20,
    text width=9em, text centered, node distance=1cm, rounded corners, minimum height=4em, font=\small]
\tikzstyle{line} = [draw, -latex']
\tikzstyle{arrow} = [draw, -latex']
\tikzstyle{cloud} = [rectangle, draw,fill=black!10, node distance=1cm, text width=8em, text badly centered,
    minimum height=2em, rounded corners, font=\footnotesize]
\begin{tikzpicture}
    [node distance = 2cm, auto
    every node/.style={node distance=3cm},
    comment/.style={rectangle, inner sep= 5pt, text width=14cm, node distance=0.25cm, font=\small},
    force/.style={rectangle, draw, fill=blue!10, inner sep=5pt, text width=4cm, text badly centered, minimum height=1.2cm, font=\bfseries\large}]
    \node [sampler_block] (coarse_sampler) { \textbf{SPDE Sampler}: Solve \eqref{eq:ml-linsys} on structured grid at coarse level $\ell+1$. };
    \node [sampler_block, below=1.75cm of coarse_sampler] (sampler) { \textbf{SPDE Sampler}: Solve \eqref{eq:ml-linsys} on structured grid at fine level $\ell$ with initial guess $\overline{\mathcal{P}} \overline{U}_{\ell+1}$. };
    \node [projection_block, right=0.8cm of coarse_sampler] (coarse_projection){\textbf{$L^2$-Projection}: \\  Compute $\btheta_{\ell+1}=\Pi_{\ell+1}\obtheta_{\ell+1}$ with coarse grid operator \eqref{eq:recursive_l2_projection}.};
    \node [projection_block, right=0.8cm of sampler] (projection){\textbf{$L^2$-Projection}:  \\ Compute $\btheta_{\ell}=\Pi_{\ell}\obtheta_{\ell}$ with fine grid operator \eqref{eq:proj_op}.};

    \node [cloud, above=0.5cm of coarse_sampler] (initial) { \footnotesize{Random Input: \\ $\xi_{\ell}(\omega_i)\sim\mathcal{N}(0,I)$}};
    \node [model_block, right=1.3cm of coarse_projection] (coarse_model){\textbf{Model evaluation}:\\ Solve \eqref{eq:ml-mixedDarcyDiscrete} on coarse level $\ell+1$ of unstructured hierarchy.};
    \node [model_block, right=1.3cm of projection] (model){\textbf{Model evaluation}: Solve \eqref{eq:ml-mixedDarcyDiscrete} on fine level $\ell$ of unstructured hierarchy.};
   
    \coordinate[below=0.75cm of coarse_model] (c);
    \node [postprocess_block, right=1cm of c] (postprocess) {\textbf{Postprocess}: Compute \\ $Q_{\ell}^{(i)}$,  $Q_{\ell+1}^{(i)}$};
    
    \coordinate[left=0.5cm of initial] (left_input);
    \coordinate[left=0.4cm of sampler] (left_sampler);
    \draw[-] (initial) -- (left_input);
    \draw[-] (left_input) -- (left_sampler);
    \path [line] (left_sampler) -- (sampler);
    \path [line] (initial) -- (coarse_sampler);
    \path [line] (coarse_sampler) -- node[anchor=east]{\footnotesize{$\overline{\mathcal{P}} \overline{U}_{\ell+1}$}}(sampler);
    \path [line] (coarse_sampler) -- node[anchor=south] {\footnotesize{$\obtheta_{\ell+1}$}} (coarse_projection);
    \path [line] (sampler) -- node[anchor=south] {\footnotesize{$\obtheta_{\ell}$}} (projection);
    
    \path [line] (coarse_projection) -- node[above] {\footnotesize{$k_{\ell+1}=$}} node[below] {\footnotesize{$\exp[\btheta_{\ell+1}]$}} (coarse_model);
    \path [line] (projection) -- node[above] {\footnotesize{$k_{\ell}=$}} node[below] {\footnotesize{$\exp[\btheta_{\ell}]$}} (model);
    
    \path [line] (model) -| node[anchor=west] {\footnotesize{$\mathbf{X}_{\ell}$}} (postprocess);
    \path [line] (coarse_model) -| node[anchor=west] {\footnotesize{$\mathbf{X}_{\ell+1}$}} (postprocess);

    \coordinate[above=1.25cm of coarse_sampler] (legend);
    \node [comment, right=1.75cm of legend] (comment-projection) {\footnotesize{$\obtheta_{\ell}$, $\obtheta_{\ell+1}$: Gaussian realizations on structured mesh \\
    $\Pi_{\ell}$, $\Pi_{\ell+1}$: $L^2$-projection operators \\ 
    $\btheta_{\ell}$, $\btheta_{\ell+1}$: Gaussian realizations on original, unstructured mesh \\
    $\mathbf{X}_{\ell}$, $\mathbf{X}_{\ell+1}$: Solutions of forward model problem } };
\end{tikzpicture}
\clearpage{}
\caption{Workflow to generate a sample, $Y_{\ell}^{(i)} = Q_{\ell}^{(i)} - Q_{\ell+1}^{(i)}$, for the MLMC estimator \eqref{eq:mlmc_estimator}, where $Q_{\ell}^{(i)}$ is the QoI on the fine level, and $Q_{\ell+1}^{(i)}$ is the QoI on the coarse level with the same random sample $\omega^{(i)}$.}
\label{fig:workflow}
\end{figure}
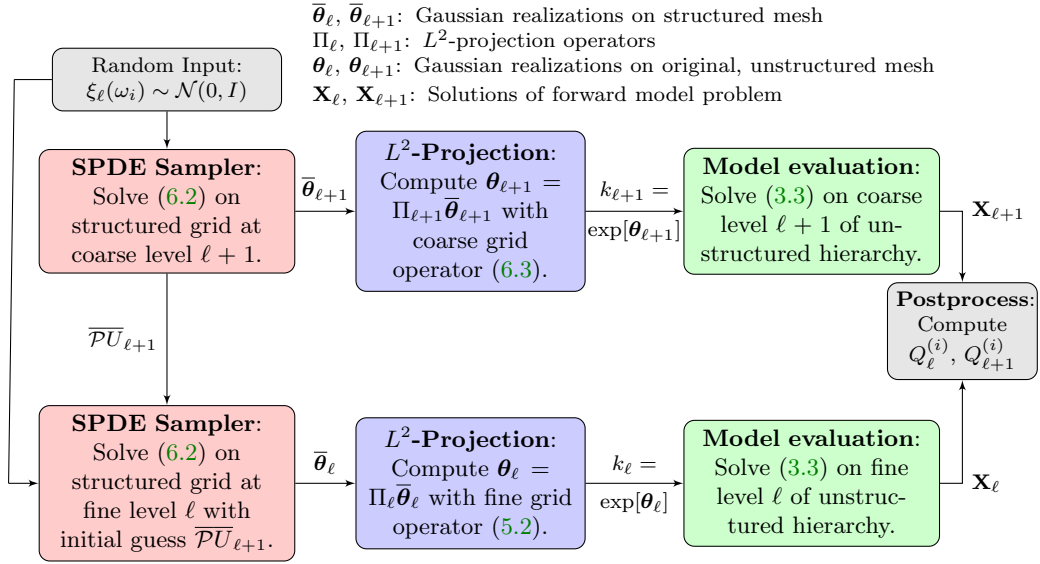

Figure~\ref{fig:workflow} illustrates the workflow to generate a MLMC sample, 
$Y_{\ell}^{(i)} = Q_{\ell}^{(i)} - Q_{\ell+1}^{(i)}$, on a particular level $\ell$
for the estimator \eqref{eq:mlmc_estimator}: First the random input coefficients 
are generated with our SPDE sampler using the non-matching mesh embedding on each level, 
then the forward model problem is solved on the    
original, unstructured grid. 

\subsection{Implementation details}\label{subsec:implementation}
We use the C++ finite element library MFEM~\cite{mfem} to assemble the
discretized problems for the sampler and forward model problem.
The hierarchy of discretizations, both structured and unstructured, 
are generated using the C++ library ParELAG~\cite{parelag}, which uses a specialized
element-based agglomeration technique to generate the algebraically constructed coarse spaces.

To formulate the right hand side of the SPDE sampler linear system 
\eqref{eq:ml_spde_sampler}, we must be able to draw a
coefficient vector of suitable independent random numbers. In our numerical
experiments,
we use Tina's Random Number Generator Library~\cite{trng}
which is a pseudo-random number generator with dedicated support for
parallel, distributed environments~\cite{bauke2007random}.

For the SPDE sampler, the saddle-point mixed linear system \eqref{eq:ml-linsys}    
is solved with the strategy described in Section~\ref{subsec:solvers} on each level, 
that is, hybridization where the reduced constrained system is solved 
with CG preconditioned with BoomerAMG~\cite{hypre}.
The forward model 
discrete saddle-point system \eqref{eq:ml-mixedDarcyDiscrete} 
is first assembled for each input realization, 
as discussed in detail in Section~\ref{subsec:ml-darcy}, 
see equations~\eqref{eq:galerkin_projection_model} and \eqref{M(k) element assembly}; then 
is solved with GMRES 
preconditioned with a block-LDU preconditioner, 
as described in Section~\ref{subsec:darcy-solver}, where
the AMG preconditioner used to apply the action of the inverse 
of the approximate Schur-complement is recomputed for each input realization.
The linear systems, both structured and unstructured,
are iteratively solved 
with an absolute stopping criterion of $10^{-12}$ and a relative stopping
criterion of $10^{-6}$.

We now investigate the performance of the hierarchical sampler under weak scaling,
i.e. when the number of mesh elements is proportional to the number of
processes, and present results for MLMC simulations for two different three-dimensional domains.
We consider a sequential, adaptive MLMC algorithm following
\cite{cliffe2011multilevel} that estimates the discretization and sampling error
from the computed samples and
chooses the optimal values for $N_\ell$ ``on the fly'' during the MLMC
simulation according to \eqref{eq:optimal_nsamples}.

The numerical experiments were executed on
\textit{quartz}, a high performance cluster at Lawrence Livermore National
Laboratory consisting of 2,688 nodes where each node has 128 GB of main memory and
36 cores operating at a clock rate of 2.1 GHz, for a total of 96768 cores.
We use the full capacity of the nodes, i.e. 36 MPI processes per node.

\subsection{Crooked pipe problem}\label{subsec:crooked_pipe_sampler}
\begin{figure}[htbp]
  \begin{subfigure}{0.3\textwidth}
    \centering
    \includegraphics[scale=.18, angle=3]{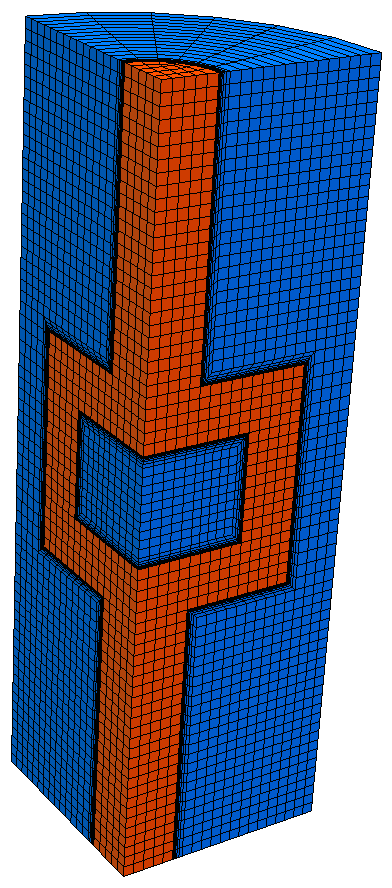}
    \caption{Crooked pipe problem}
  \end{subfigure}
  \begin{subfigure}{0.3\textwidth}
    \centering
    \includegraphics[scale=.27]{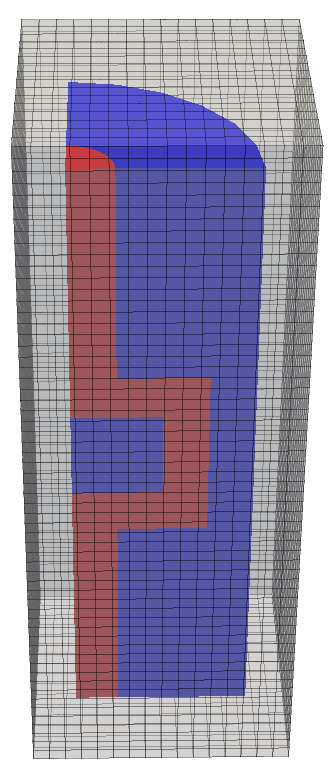}
    \caption{Regular enlarged mesh embedding}
  \end{subfigure}
  \begin{subfigure}{0.35\textwidth}
    \centering
    \includegraphics[scale=.1]{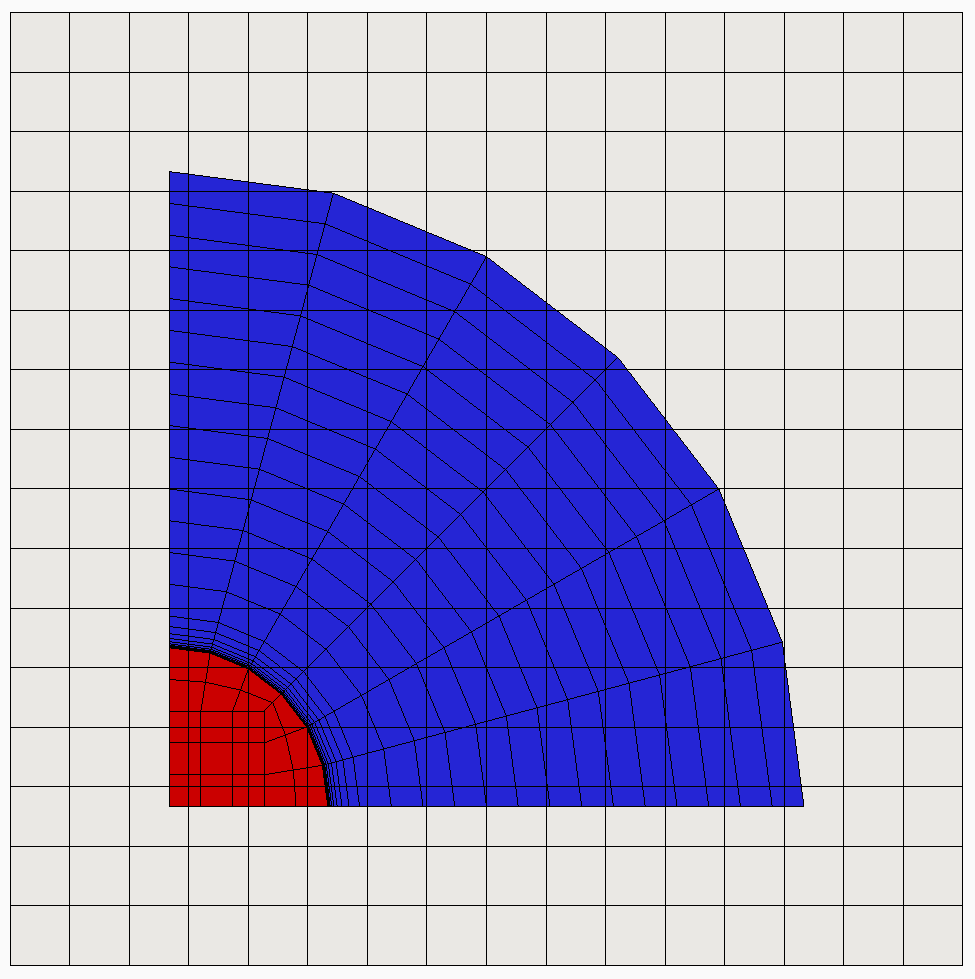}
    \caption{Slice view of non-matching mesh structure}
    \centering
    \includegraphics[scale=.23]{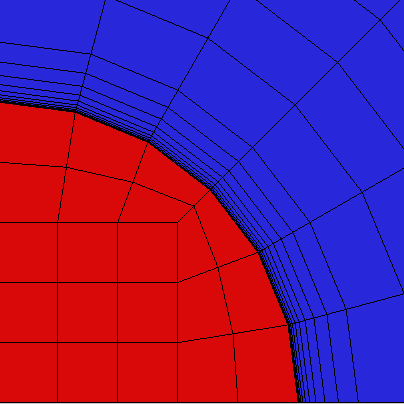} 
    \caption{Zoom view of highly stretched elements}
  \end{subfigure}
  \caption{ The initial crooked pipe mesh $D$, a quarter cylinder 
      shape with radius equal to 2 and height equal to 7, featuring highly
      stretched elements that are used to capture the boundary layer at the
      interface between two material subdomains, with $14370$ hexahedral elements and
      the larger, regular bounding box
        $\oD=(0,3)\times(0,3)\times(0,8)$ with $15360$ hexahedral elements.} 
  \label{fig:crooked_pipe_mesh}
\end{figure}

\begin{figure}[htbp]
  \begin{subfigure}{0.24\textwidth}
    \centering
    \includegraphics[scale=.25]{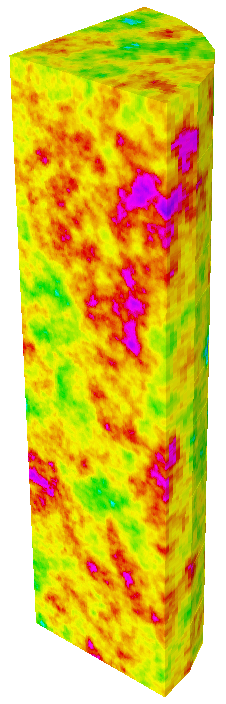}
    \caption{Level $\ell=0$}
  \end{subfigure}
  \begin{subfigure}{0.24\textwidth}
    \centering
    \includegraphics[scale=.25]{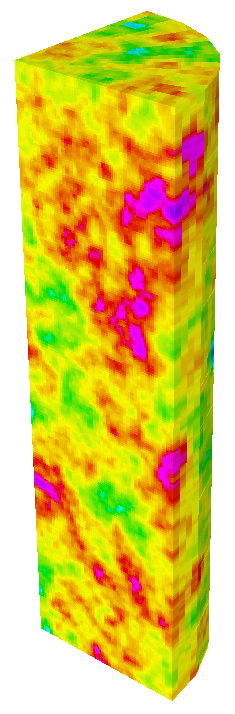}
    \caption{Level $\ell=1$}
  \end{subfigure}
  \begin{subfigure}{0.24\textwidth}
    \centering
    \includegraphics[scale=.25]{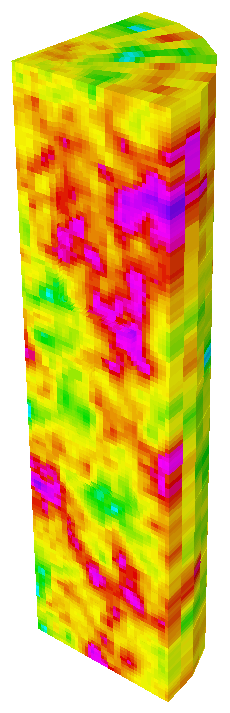}
    \caption{Level $\ell=2$}
  \end{subfigure}
  \begin{subfigure}{0.24\textwidth}
    \centering
    \includegraphics[scale=.25]{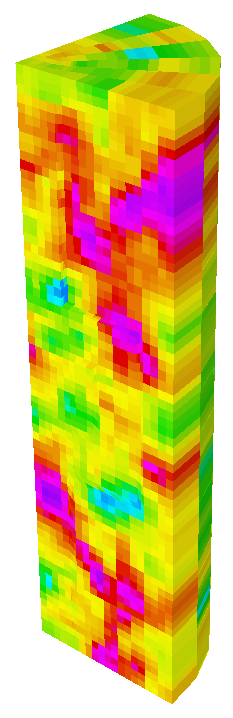}
    \caption{Level $\ell=3$}
  \end{subfigure}
  \caption{A realization of the Gaussian random field on the crooked pipe domain obtained by using our
     hierarchical sampling technique for 4 levels with \Mat covariance with correlation length $\lambda=0.3$. }
  \label{fig:crooked_pipe_realizations}
\end{figure}

We first consider a cylindrical ``butterfly''-type grid, 
with highly 
stretched elements that are used to capture the boundary layer at the 
interface between two material subdomains. 
The domain $D$, a quarter cylinder with radius equal 2 and height equal 7,
is embedded in the regular grid given by $\oD=(0,3)\times(0,3)\times(0,8).$
Each mesh is uniformly refined several times to build the hierarchy of levels.
Figure~\ref{fig:crooked_pipe_mesh} shows the initial mesh for the crooked pipe
problem and the enlarged, regular domain $\oD$.
  The parameters of the random field are variance $\sigma^2 = 1$ and correlation
  length $\lambda = 0.3$. A sequence of computed Gaussian field realizations on
  different levels is shown in Figure~\ref{fig:crooked_pipe_realizations}.

\begin{figure}[htbp] \centering
  \begin{subfigure}{0.6\textwidth}
    \scalebox{0.9}{
    \begin{tabular}{lrrr}
    \toprule
    Processes&  144&  1152&  9216  \\
    \midrule
    Construct $\Pi_0$& 27.1485 & 29.7894 & 42.0257\\ 
    Preconditioner Set-up& 2.7156& 2.8003 & 2.9628\\
    Solve \eqref{eq:linsys} for $\obtheta_0$ & 140.6940 & 157.1750 & 191.7121\\ 
    Compute $\btheta_0=\Pi_0\obtheta_0$ & 1.3203 & 1.2112 & 0.9331\\
    \bottomrule
    \end{tabular} 
    } 
    \caption{Computational time (secs) of generating 100 Gaussian realizations
    on the fine level using SPDE sampler.}
    \label{subfig:crooked_pipe_sampler_table}
  \end{subfigure}
  \begin{subfigure}{0.38\textwidth}
    \includegraphics[scale=.35]{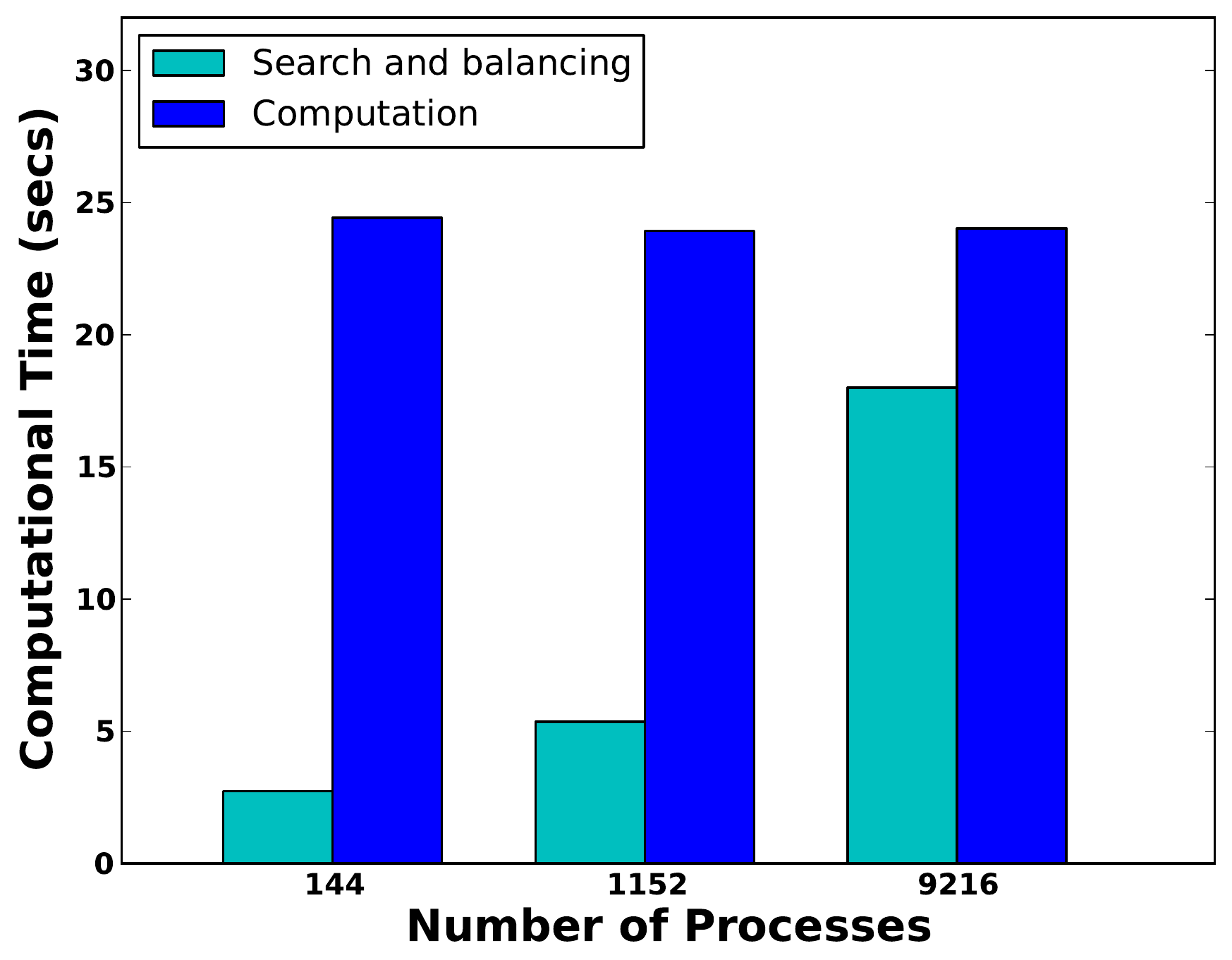}
    \caption{Overview of the search and computation times for the assembly of the $L^2$-projection.}
    \label{subfig:crooked_pipe_l2proj}
  \end{subfigure}
  \begin{subfigure}{0.99\textwidth}
  \begin{center}
  \begin{tabular}{lrrr}
  \toprule
  Processes&  144&  1152&  9216  \\
  \midrule
  Element bounding volumes generation & 0.1783 & 0.1788 & 0.1783\\
  BVH comparison & 0.4726 & 0.6303 & 2.5610 \\
  Load balancing & 0.6702 & 1.2605 & 1.385 \\
  Matching and rebalancing & 1.4066 & 3.7860  & 13.8728 \\
  \hdashline
  Computation: intersection and quadrature & 24.4207 & 23.9338 & 24.0251\\
  \hline
   Total  & 27.1485 & 29.7894  & 42.0257\\
  \bottomrule
  \end{tabular}\end{center}
    \caption{
    Detailed computational times for the assembly of the $L^2$-projection. Listed above the dashed-line
    are the different components of the \emph{search and balancing} phase.
  }
  \label{subfig:crooked_pipe_sampler_details_table}
  \end{subfigure}
  
  \caption{The computational cost for the crooked pipe mesh 
    of generating 100 Gaussian realizations on the fine level 
    under weak scaling with approximately $5.1\cdot 10^4$ stochastic degrees of 
    freedom per process is shown in (a).
  The weak scalability of the $L^2$-projection operator assembly 
  where the input is approximately $1.1\cdot 10^5$ elements per process and
  the average number of intersections per process is $3.0\cdot 10^5$ is demonstrated 
  in (b) and (c). 
  The computational time is divided between a \emph{search and
  balancing} phase and a phase of \emph{computation} of the transfer operator in (b)
  with the detailed timing results exhibited in (c).} 
  \label{fig:crooked_pipe_sampler_setup}
\end{figure}

We examine the performance of the SPDE sampler under weak scaling for the crooked pipe domain
in Figure~\ref{fig:crooked_pipe_sampler_setup}. The computational
time to construct the SPDE sampler and generate 100 Gaussian realizations on the fine
level with $5.1\cdot 10^4$ stochastic degrees of freedom per
process is examined in Figure~\ref{subfig:crooked_pipe_sampler_table}.
In our approach the
construction time of the $L^2$-projection operator takes about $16-18\%$ of the
total computational time to generate 100 Gaussian realizations (Figure~\ref{subfig:crooked_pipe_sampler_table}).
Weak scaling computational times to construct the
$L^2$-projection operator with $1.1\cdot 10^5$ input elements per process on the fine level
are reported in Figure~\ref{subfig:crooked_pipe_l2proj} for both the 
\emph{search and balancing} and \emph{computation} phase, and the 
detailed computational times are shown in Figure~\ref{subfig:crooked_pipe_sampler_details_table} 
where \emph{search and balancing} includes all measurements except the computation.
The timing labels are described in Section~\ref{sub-sec:mapping}.

The \emph{search and balancing} phase finds intersecting elements and 
redistributes elements to ensure load balancing of the
computation. The search is a global process which requires communication and
synchronization which leads to additional overhead when adding more processes.
As described in Section~\ref{sub-sec:mapping}, the search algorithm is based
on the construction of an octree. The cost of this step
can be tuned by the user by changing tree construction parameters such
as maximum tree-depth and maximum number of elements per tree-node, however the computation has a
lower bound computational time complexity $\Omega(n \log_8 n)$, where $n$ is the number of input
elements. From this lower bound we can expect that when increasing the size of
the input by one order of magnitude we loose approximately $50\%$ weak-scaling
efficiency when searching for intersecting element pairs (Figure~\ref{subfig:crooked_pipe_sampler_details_table}).
These limitations are represented in the measurements. With loss of generality
the intersection detection could be further optimized for dealing with
Cartesian grids where usually the (implicit) grid information is globally
known~\cite{cavoretto2017opencl}. 
In the \emph{computation} phase we calculate polyhedral intersections
and perform numerical quadrature. Here we also determine the exact size of the
output which cannot be predicted accurately in advance, hence we encounter
unavoidable slight computational imbalances.

\begin{figure}[htbp] \centering
  \begin{subfigure}{0.49\textwidth}
    \includegraphics[scale=.36]{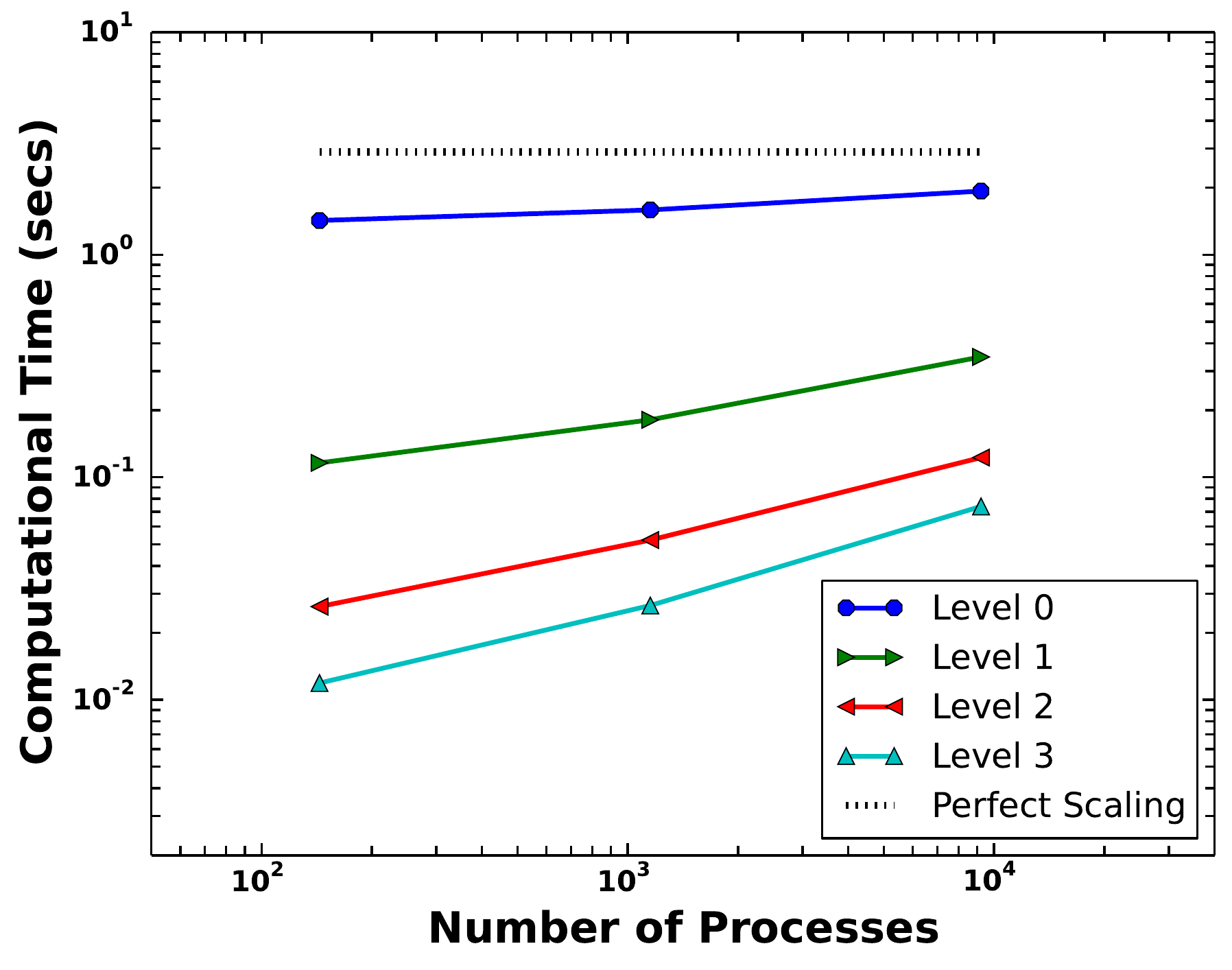}
    \label{subfig:crooked_pipe_sampler_weak_scaling_solve_time}
  \end{subfigure}
  \begin{subfigure}{0.49\textwidth}
    \includegraphics[scale=.36]{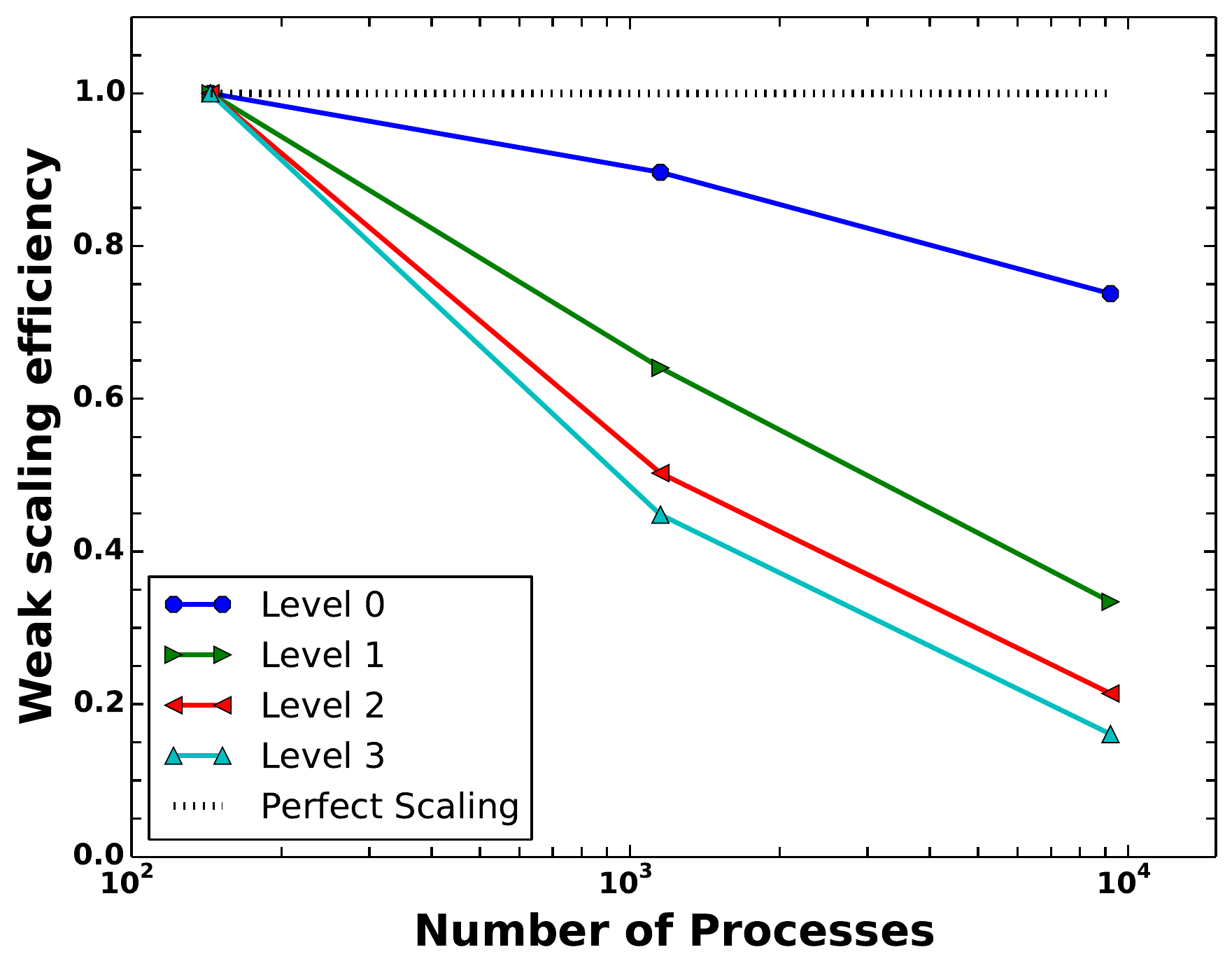}
    \label{subfig:crooked_pipe_sampler_efficiency}
  \end{subfigure}
  \caption{ Weak scalability of the linear solve time and $L^2$-projection operator 
    applied for the hierarchical SPDE sampler on the crooked pipe
    domain where the size of the stochastic dimension per
    process, approximately $5.1\cdot 10^4$ stochastic degrees of
    freedom on the fine level, is fixed. The size of the stochastic
    dimension of the finest level ranges from $7.4\cdot 10^6$ to $4.7\cdot
    10^8$, and the number of MPI processes ranges from 144
    to 9216. We observe approximately 74\% parallel efficiency 
    for the fine level with 9216 processes to compute a realization using the sampling method.}
  \label{fig:crooked_pipe_sampler_weak_scaling}
\end{figure}

Figure~\ref{fig:crooked_pipe_sampler_weak_scaling} shows the average time --- 
for a number of MPI processes ranging from 144
to 9216 --- to compute
a realization on $\oD$ and project the solution to the crooked pipe domain, $D$, with
approximately $5.1\cdot 10^4$ stochastic degrees of freedom per
process on the fine level. 
We observe approximately 74\% parallel efficiency
for the fine level with 9216 processes to compute a realization using the sampling method
with the scalable hybridization multigrid preconditioning strategy.

Now we consider a MLMC simulation with the forward model given by 
\eqref{eq:darcy} with boundary conditions given by
    \begin{equation*}
        \begin{cases}
            -p&=1 \quad \text{on} \ \Gamma_{in}, \\
            -p&=0 \quad \text{on} \ \Gamma_{out}, \\
            \bq \cdot {\bf n} &= 0\quad  \text{on} \ \Gamma_{s} := \partial D \setminus \left(\Gamma_{in} \cup \Gamma_{out}\right), \\
        \end{cases}
    \end{equation*}
    where $\Gamma_{in}$ is the boundary along the plane at $z=0.5$ 
    and $\Gamma_{out}$ is the boundary along the plane at $z=7.5$.

\begin{figure}[htbp]
  \begin{subfigure}{0.5\textwidth}
    \centering
    \includegraphics[scale=.34]{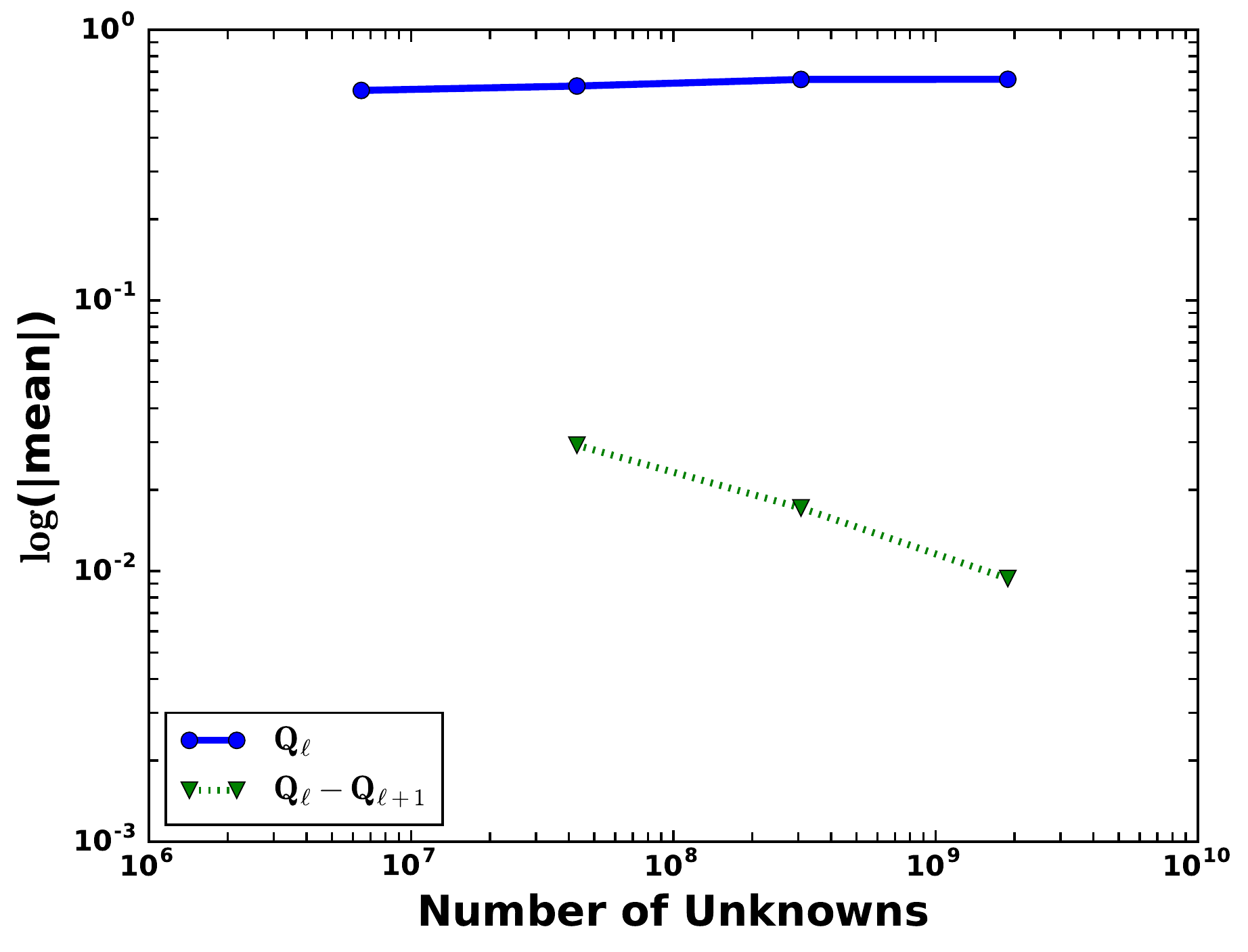}
   \caption{Multilevel estimator}
   \label{subfig:crooked_pipe_mean}
  \end{subfigure}
  \begin{subfigure}{0.5\textwidth}
    \centering
    \includegraphics[scale=.34]{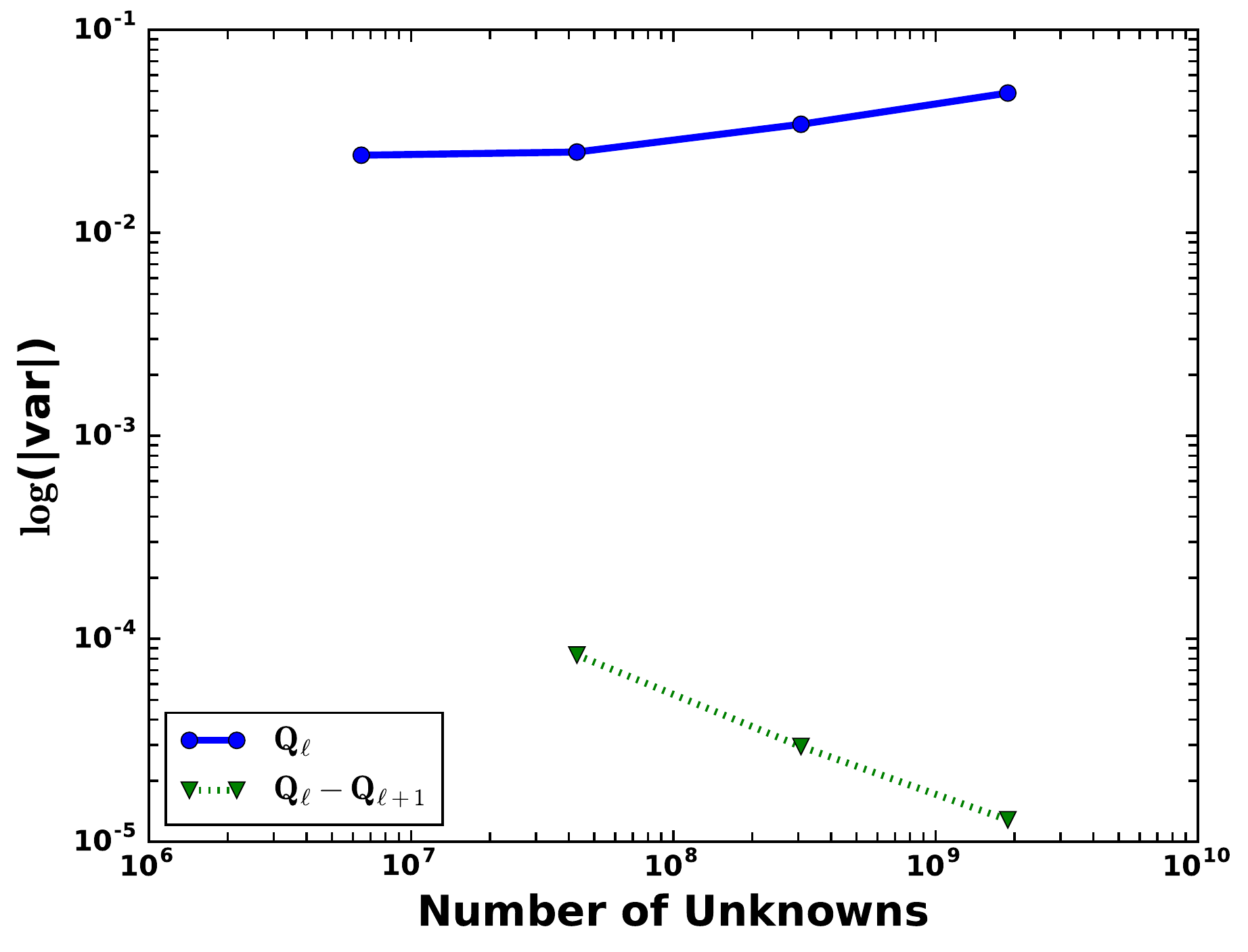}
   \caption{Variance reduction}
   \label{subfig:crooked_pipe_var}
  \end{subfigure}
  \begin{subfigure}{0.49\textwidth}
    \centering
    \includegraphics[scale=.34]{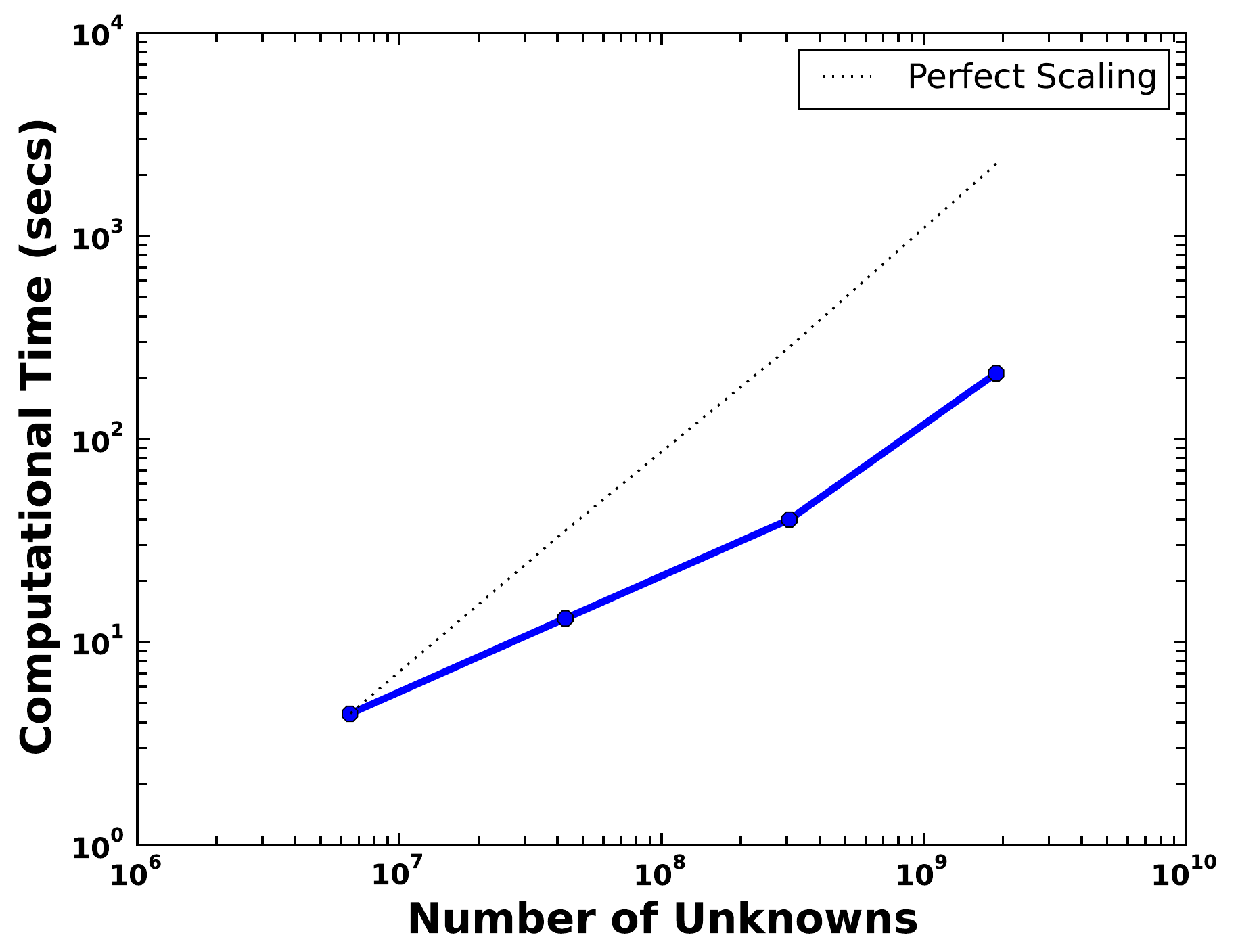}
   \caption{Average sample time}
   \label{subfig:crooked_pipe_time}
  \end{subfigure}
  \begin{subfigure}{0.49\textwidth}
    \centering
    \includegraphics[scale=.33]{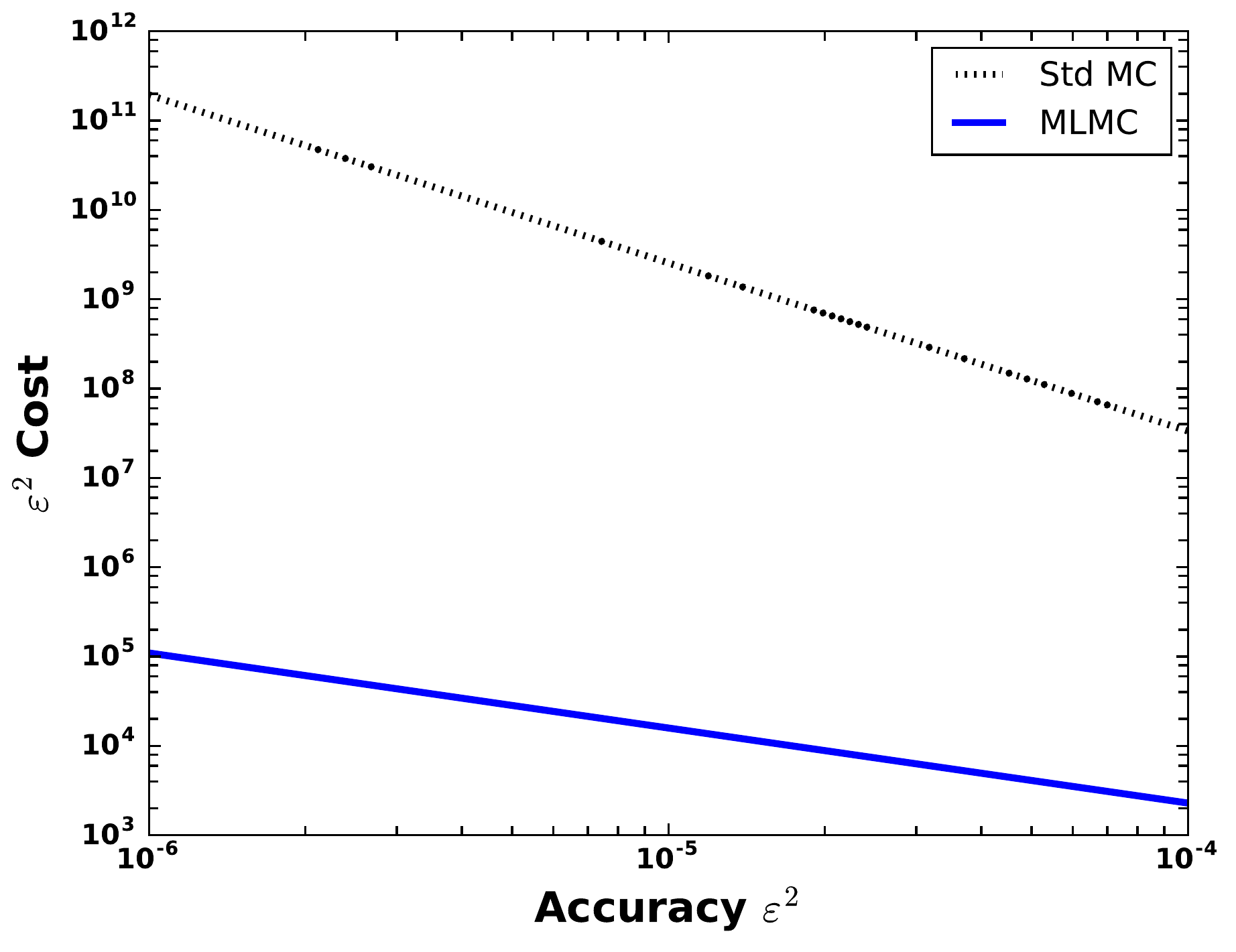}
   \caption{Predicted Asymptotic Order of Cost}
   \label{subfig:crooked_pipe_theoretical}
  \end{subfigure}
 
  \begin{subfigure}{0.99\textwidth}
  \begin{center}
  \scalebox{0.9}{
  \begin{tabular}{lrrrr}
  \toprule
  Level &  0 (fine)&  1&  2& 3 (coarse)  \\
  \midrule
  Number of samples & 42 & 99 & 234 & 7678 \\
  \hline 
  Realize input coefficients (SPDE sampler) & 95.6791 & 46.5105 & 45.9743 & 567.5747\\
  Assemble systems \eqref{eq:ml-mixedDarcyDiscrete} on levels $\ell$ and $\ell+1$& 315.2517 & 123.2833 & 20.1287 & 86.9110  \\
  Build AMG preconditioner for $\tilde{S}_{\ell}$ and $\tilde{S}_{\ell+1}$ in \eqref{eq:ldu} & 4.1919 & 2.5160 & 0.8775 & 10.3692 \\
  Solve \eqref{eq:ml-mixedDarcyDiscrete} on levels $\ell$ and $\ell+1$& 2516.1402& 3557.3076 & 4172.4751 & 19903.6026\\
  \hline
   Total time & 2931.2629 & 3729.6174  & 4239.4555 & 20568.4575\\
  \bottomrule
  \end{tabular} }
  \end{center}
    \caption{
    Detailed computational times (secs) for each level of the MLMC simulation. 
  }
  \label{subfig:crooked_pipe_mlmc_table}
  \end{subfigure}
  \caption{MLMC results for crooked pipe problem when estimating the
      effective permeability, where the target MSE is
      $\varepsilon^2=6.25\cdot 10^{-6}$ using $9216$ processes with 
      approximately $1.9 \cdot 10^9$ velocity and pressure degrees of freedom on the
      fine level. The expected value of the MC estimator of $Q_{\ell}$ and 
      the correction are shown in
      (a) and the variance reduction of the multilevel method 
      is shown in (b) where the variance of the multilevel correction
      term is significantly smaller than the variance of the MC estimator of $Q_{\ell}$.
      In (c), the average sampling
      time to generate a MLMC sample $Y_{\ell}^{(i)}=Q_{\ell}^{(i)}-Q_{\ell+1}^{(i)}$ 
      for each level versus the number of
      unknowns on level $\ell$ is plotted. Plot (d) compares the theoretical asymptotic order of cost to
      achieve a MSE of $\varepsilon^2$ for this problem formulation for the standard
      MC and MLMC estimators. The MLMC method leads to a
      significant improvement over the standard MC method. The time spent on 
      each level of the MLMC hierarchy is shown in (e).}\label{fig:crooked_pipe_mlmc}
\end{figure}

The random input coefficient is assumed to be log-normal with \Mat covariance
    (equivalent to an exponential covariance since $d=3$) with $\sigma^2 = 1$ and correlation
    length $\lambda = 0.3$. The hierarchical SPDE sampler with mesh embedding
    is used to generate the input realizations.
    The quantity of interest is the expected value of the
effective permeability, that is the flux through the ``outflow''
part of the boundary, defined as
\begin{equation}
    k_{eff}(\omega) = \frac{1}{\abs{\Gamma_{out}}} \int_{\Gamma_{out}} \bq 
    (\cdot,\omega)\cdot{\bf n} \, dS.
\label{eq:eff_perm}
\end{equation}
In Figure~\ref{fig:crooked_pipe_mlmc}, we present standard MLMC results for a $4$-level method
using $9216$ processes with approximately $1.9 \cdot 10^9$ velocity and pressure degrees of freedom on the fine level 
and target MSE $\varepsilon^2=6.25\cdot 10^{-6}$.
Figure~\ref{subfig:crooked_pipe_mean} displays the multilevel MC estimator, where
the blue solid line represents the
expectation at each level, $\E[Q_{\ell}]$, and the green dashed line  
represents the expectation of the correction, $\E[Q_{\ell}-Q_{\ell+1}]$.
Figure~\ref{subfig:crooked_pipe_var} illustrates the multilevel variance reduction
where the blue solid line shows the variance of the estimator for a particular
level, whereas the green dashed line shows the variance of the correction 
for each level. This plot demonstrates the effectiveness of the MLMC method:
as the number of unknowns (i.e. the spatial resolution) increases, the variance of 
the correction is reduced.
The average sampling time to generate the required Gaussian field realizations
and solve the forward model for each level is shown in
Figure~\ref{subfig:crooked_pipe_time}. This plot indicates near optimal scaling of the MLMC
method with the proposed hierarchical sampler.
Figure~\ref{subfig:crooked_pipe_theoretical} compares the predicted $\varepsilon^2$-cost of
the standard MC and MLMC estimators using Theorem 2.3 from~\cite{teckentrup2013further} with
numerically observed constants estimated from
Figures~\ref{subfig:crooked_pipe_mean}-\ref{subfig:crooked_pipe_time}.      
The detailed computational time spent on each level is shown in Figure~\ref{subfig:crooked_pipe_mlmc_table},
where the majority of the time is spent generating samples on the coarsest level. 

Next we consider the performance of a MLMC simulation for the crooked pipe problem
under weak scaling, where the number of velocity and pressure
degrees of freedom per process is fixed. Table~\ref{tab:crooked_pipe_weak_table} displays
the results of three MLMC simulations with increasing spatial resolution along
with the processor count, while the desired MSE is decreasing. The
tolerance for the sampling error is chosen to balance the estimated
discretization error so more samples are necessary as a finer spatial resolution is used,
yet the majority of the samples are computed on the coarsest level.
\begin{table}[tbhp]
  \caption{Weak scaling of the crooked pipe problem for adaptive MLMC simulations
  where the number of degrees of freedom per
  process is kept approximately fixed with approximately 
  $2.1\cdot10^5$ pressure/velocity degrees of freedom per process
  on the fine level.
  In each row, the number of MPI processes are listed along with the
  number of velocity and pressure degrees of freedom on the
  fine level ($\ell=0$), the desired MSE,
  the total computational wall time to run the MLMC simulation,
  the number of samples computed on the fine level ($N_0$), and the total
  number of computed samples on all levels for a 4 level method.}
  \label{tab:crooked_pipe_weak_table}
  \begin{center}
  \begin{tabular}{lllllr}
    \toprule
    Processes&  DOF $(\ell=0)$& Target $\varepsilon^2$ & Wall Time (s)&  $N_0$&  Total Samples \\
    \midrule
    144&  2.96$\cdot 10^7$&  $1.00\cdot 10^{-4}$    &  2.93$\cdot10^2$ &  10&  952 \\
    1152&  2.36$\cdot 10^8$& $2.50\cdot 10^{-5}$  &  1.65$\cdot10^3$ &  21&  3634 \\
    9216&  1.88$\cdot 10^9$& $6.25\cdot 10^{-6}$ &  3.15$\cdot10^4$&  42&  8053 \\
    \bottomrule
  \end{tabular}
  \end{center}
  \vspace{1ex}
\end{table}

\subsection{SPE10 problem}\label{subsec:spe10_sampler}
Next we consider the domain from Model 2 of the tenth SPE comparative solution project
(SPE10)~\cite{spe_data}, a challenging benchmark for reservoir simulation codes.
The domain is a 3D box with dimension $1200 \times 2220 \times 170$(ft) meshed with
hexahedral elements. The mesh is embedded in a bounding box with dimension
$1600 \times 2420 \times 240$(ft). 

\begin{figure}[htbp] \centering
  \begin{subfigure}{0.6\textwidth}
    \begin{center}
    \scalebox{0.9} { \begin{tabular}{lrrr}
    \toprule
    Processes&  36&  288&  2304  \\
    \midrule
    Construct $\Pi_0$& 6.9641 &  7.2735 & 7.6082\\ 
    Preconditioner Set-up& 2.0416 & 2.1059 & 2.2118\\
    Solve \eqref{eq:linsys} for $\obtheta_0$& 91.2868 & 102.9320 &135.1250\\ 
    Compute $\btheta_0=\Pi_0\obtheta_0$& 0.5786 &0.3380  &0.4868\\
    \bottomrule
    \end{tabular} } \end{center}
    \caption{Computational time (secs) of computing 100 Gaussian
    realizations on the fine level using SPDE sampler. }
    \label{subfig:SPE10_sampler_table}
  \end{subfigure}
  \begin{subfigure}{0.38\textwidth}
    \includegraphics[scale=.35]{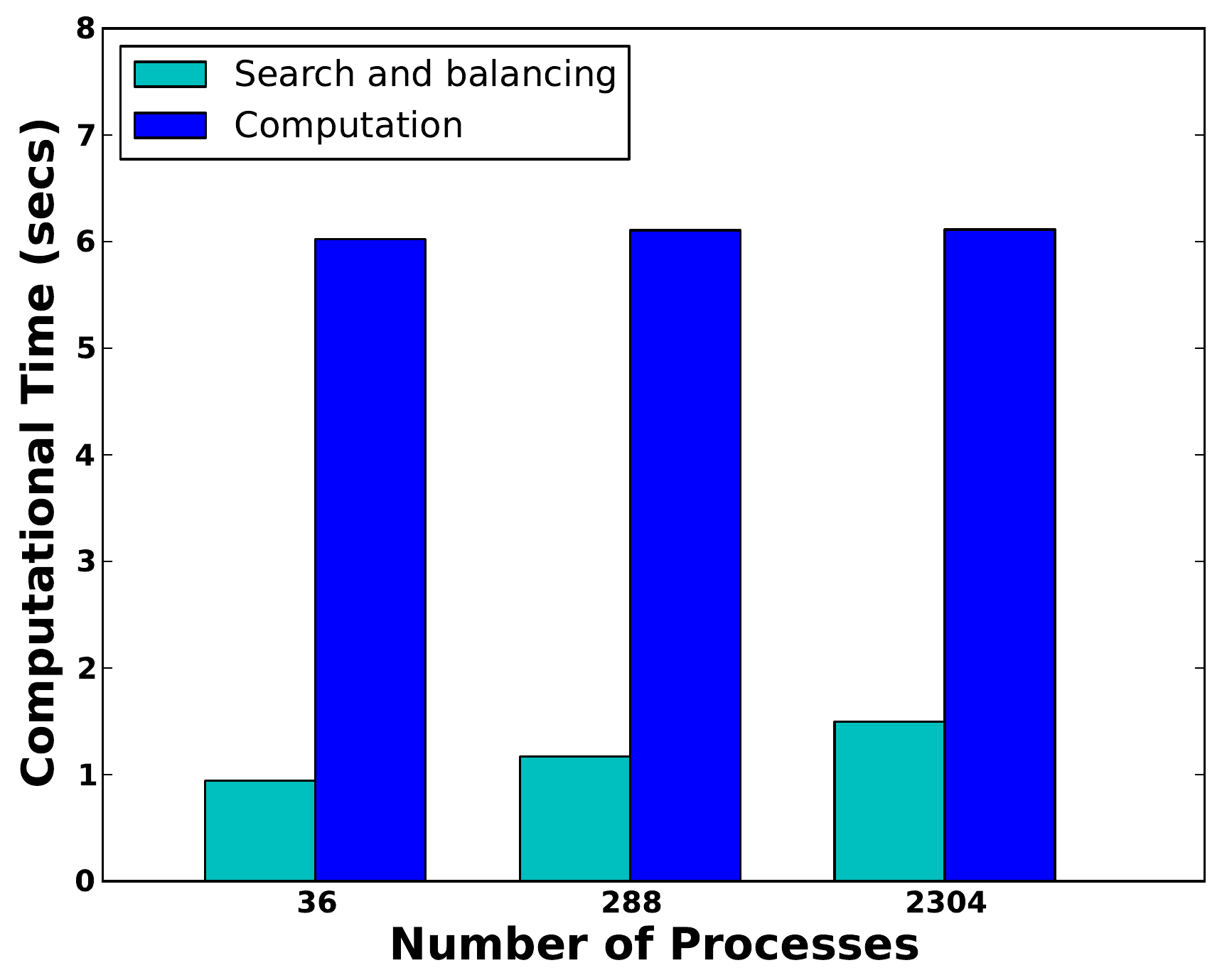}
    \caption{Overview of the search and computation times for the assembly of the $L^2$-projection.}
    \label{subfig:SPE10_l2proj}
  \end{subfigure}
  \begin{subfigure}{0.99\textwidth}
  \begin{center}\begin{tabular}{lrrr}
  \toprule
  Processes&  36&  288&  2304 \\
  \midrule

  Element bounding volumes generation & 0.0702& 0.0864& 0.0873 \\
  BVH comparison & 0.3738 & 0.3203 & 0.4937\\
  Load balancing & 0.3813 &  0.3704 & 0.3476\\
  Matching and rebalancing & 0.1174 &0.3919 & 0.5664\\
  \hdashline
  Computation: intersection and quadrature & 6.0213& 6.1045&  6.1132\\
  \hline
  Total & 6.9641 & 7.2735 & 7.6082\\
  \bottomrule
  \end{tabular}\end{center}
  \caption{
    Detailed computational times for the assembly of the $L^2$-projection. Listed above the dashed-line
    are the different components of the \emph{search and balancing} phase.}
  \label{subfig:SPE10_sampler_details_table}
  \end{subfigure}
  \caption{ The computational cost of generating 100 Gaussian realizations 
    for $D=1200\times2200\times170$ on the fine level
    under weak scaling with approximately $6.2\cdot 10^4$ stochastic degrees of
    freedom per process is shown in (a). Weak scalability of the $L^2$-projection operator assembly
   with approximately $4.8\cdot 10^4$ input elements per process
  where the average number of intersections per process is $1.5\cdot 10^4$ is shown in (b) and (c). 
  In (b), the computational time is divided between a \emph{search and
balancing} phase and a phase of \emph{computation} of the transfer operator, where the detailed 
timing results are presented in (c).}
  \label{fig:SPE10_sampler_setup}
\end{figure}

We first examine the performance of the SPDE sampler under weak scaling and 
assume the random field has $\sigma^2=1$ and correlation length
$\lambda=50$(ft).
In Figure~\ref{fig:SPE10_sampler_setup}, we examine the set-up costs associated with 
the SPDE sampler and the computational time to generate 100 samples on the fine 
level with approximately $6.2\cdot 10^4$ stochastic degrees of freedom per process. 
In our approach the
construction time of the $L^2$-projection operator takes about $5-6\%$ of the
total computational time to generate 100 Gaussian realizations 
(Figure~\ref{subfig:SPE10_sampler_table}). 
Using the labels described in Section~\ref{sub-sec:mapping}
for the computational components, the weak scaling of the construction of the
$L^2$-projection operator is shown in Figures~\ref{subfig:SPE10_l2proj} 
and~\ref{subfig:SPE10_sampler_details_table}, 
which displays a scaling behavior consistent with
the results and discussion covered in
Section~\ref{subsec:crooked_pipe_sampler} for the crooked pipe problem.

\begin{figure}[htbp] \centering
  \begin{subfigure}{0.49\textwidth}
    \includegraphics[scale=.36]{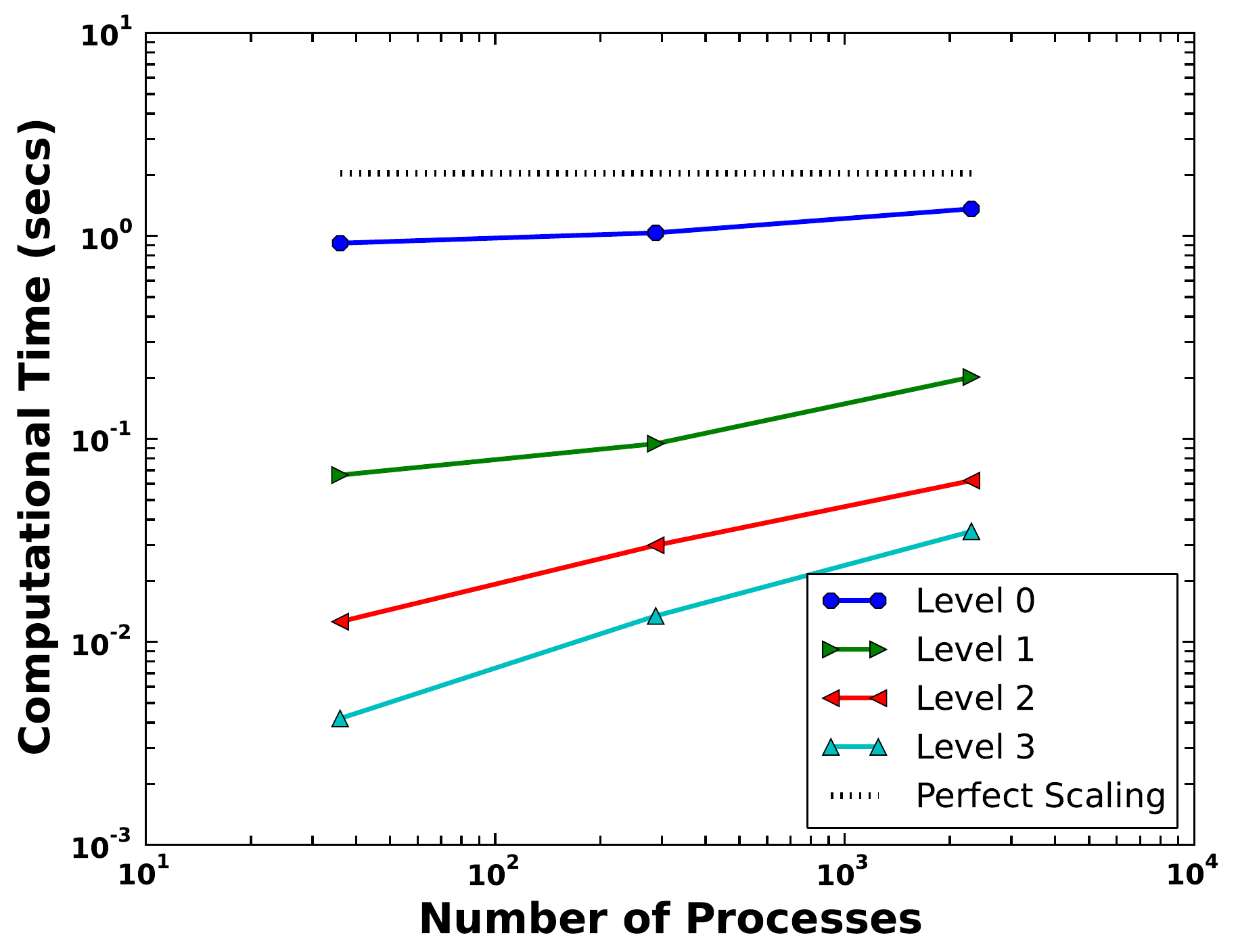}
  \label{subfig:SPE10_sampler_weak_scaling}
  \end{subfigure}  
  \begin{subfigure}{0.49\textwidth}
    \includegraphics[scale=.36]{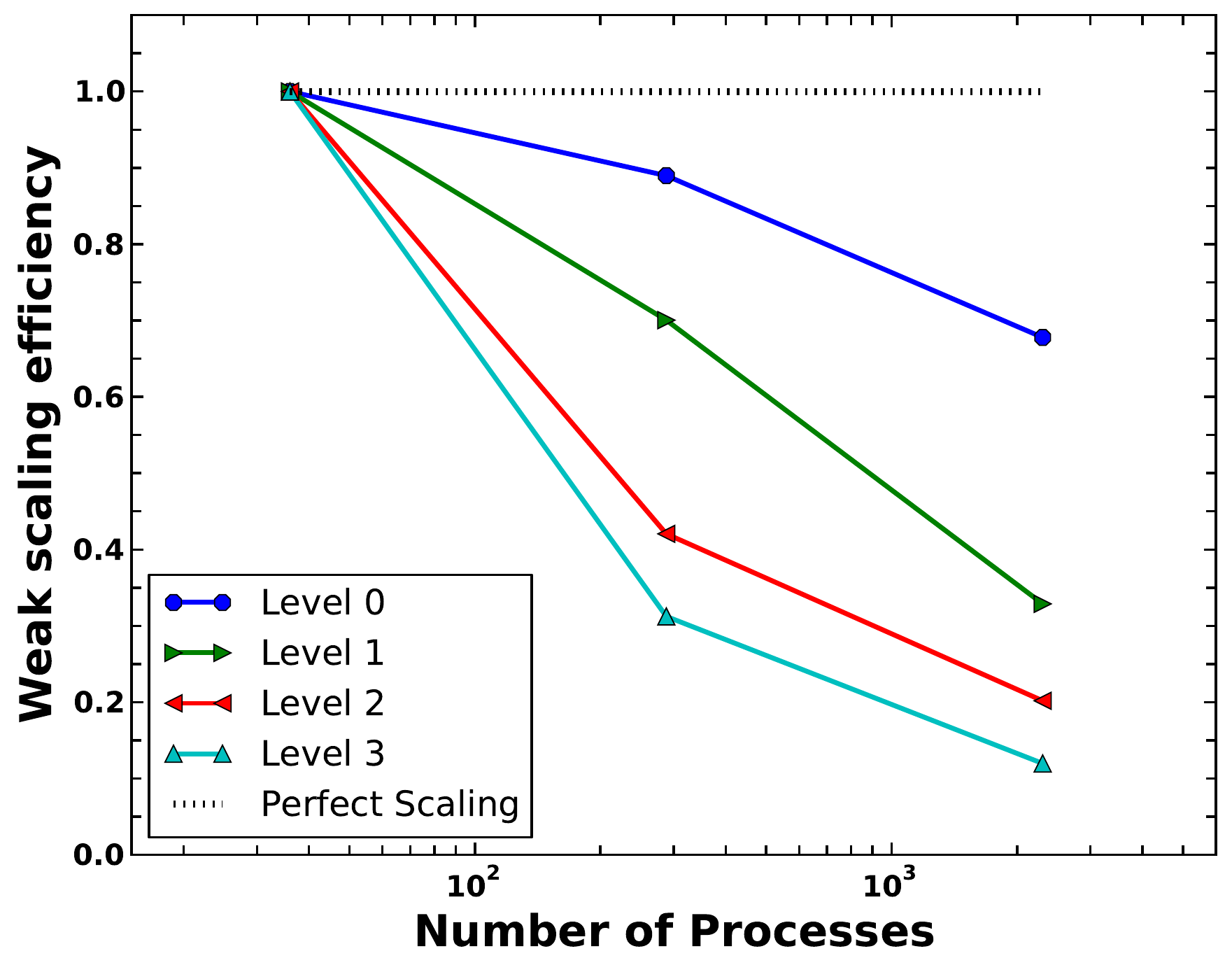}
  \label{subfig:SPE10_sampler_efficiency}
  \end{subfigure}  
  \caption{ Weak scalability of the hierarchical SPDE sampler on the domain with dimension
    $1200 \times 2220 \times 170$(ft) embedded in a box with dimension $1600 \times 2420 \times 240$(ft).
    The size of the stochastic dimension per
    process is fixed with approximately $3.3\cdot 10^4$ stochastic degrees of freedom per process
    on the fine level. The size of the stochastic dimension of the finest level
    ranges from $1.2\cdot 10^6$ to $7.5\cdot 10^7$ and the number of MPI processes ranges from 36 
    to 2304. The average time to compute a realization using the sampling method
    exhibits 68\% parallel efficiency for 2304 processes on the fine level.}
  \label{fig:SPE10_sampler_weak_scaling}
\end{figure}

Figure~\ref{fig:SPE10_sampler_weak_scaling} demonstrates the linear solver and $L^2$-projector application
performance of the
sampler under weak scaling with approximately $6.2\cdot 10^4$ stochastic degrees of freedom per process
on the fine level where the number of MPI processes ranges from 36 
to 2304. The average time to compute a realization using our sampling method
exhibits 68\% parallel efficiency for 2304 processes on the fine level.

\begin{figure}[htbp]
  \begin{subfigure}{0.49\textwidth}
    \centering
    \includegraphics[scale=.2]{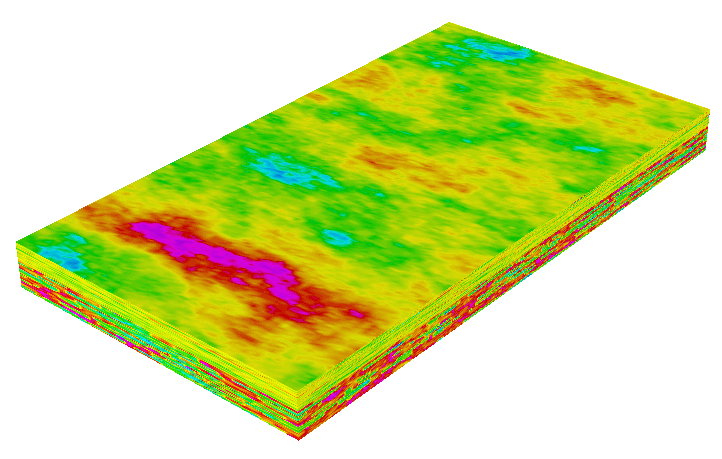}
   \caption{x/y-component}
   \label{subfig:spe10_x}
  \end{subfigure}
  \begin{subfigure}{0.49\textwidth}
    \centering
    \includegraphics[scale=.2]{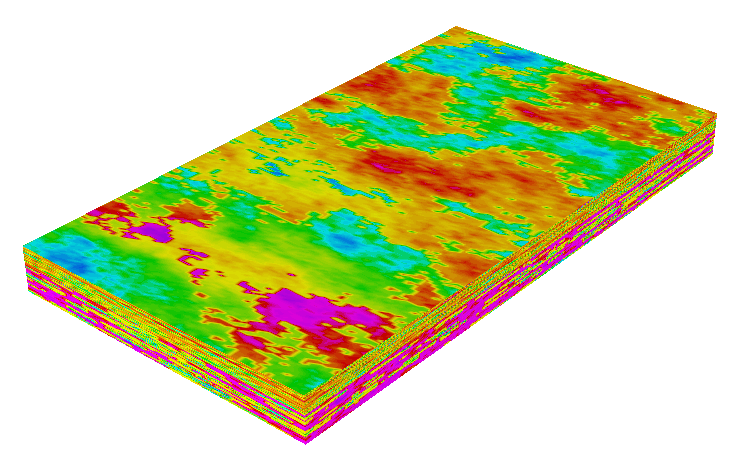}
   \caption{z component}
   \label{subfig:spe10_z}
  \end{subfigure}
  \caption{Logarithmic plots of the absolute permeability coefficient from the
      SPE10 dataset which represents the mean of the log-normal random field used to model the
    random permeability field.}\label{fig:spe10_perm}
\end{figure}

\begin{figure}[htbp]
  \begin{subfigure}{0.49\textwidth}
    \centering
    \includegraphics[scale=.34]{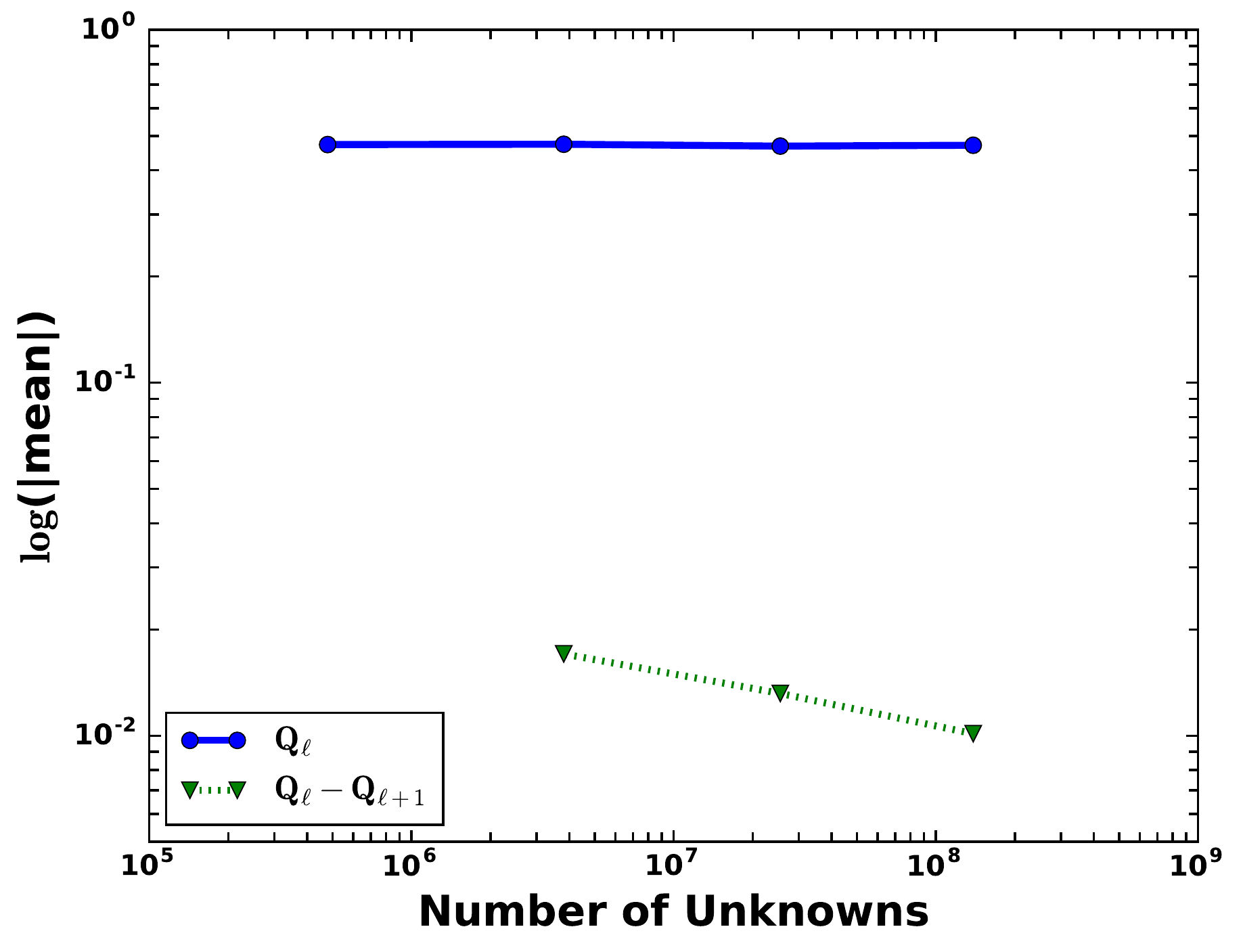}
   \caption{Multilevel estimator}
   \label{subfig:SPE10_mean}
  \end{subfigure}
  \begin{subfigure}{0.49\textwidth}
    \centering
    \includegraphics[scale=.34]{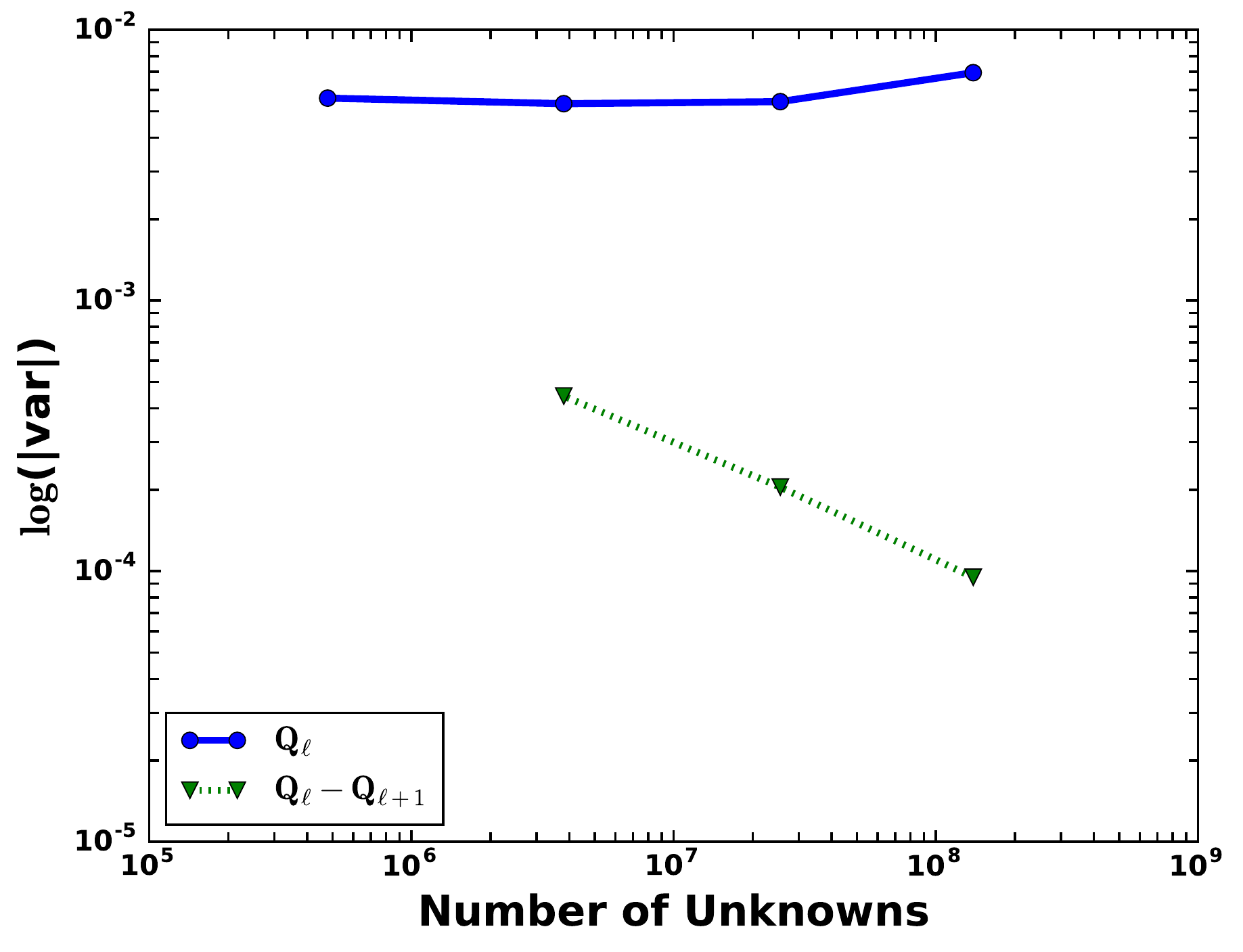}
   \caption{Variance reduction}
   \label{subfig:SPE10_var}
  \end{subfigure}
  \begin{subfigure}{0.49\textwidth}
    \centering
    \includegraphics[scale=.34]{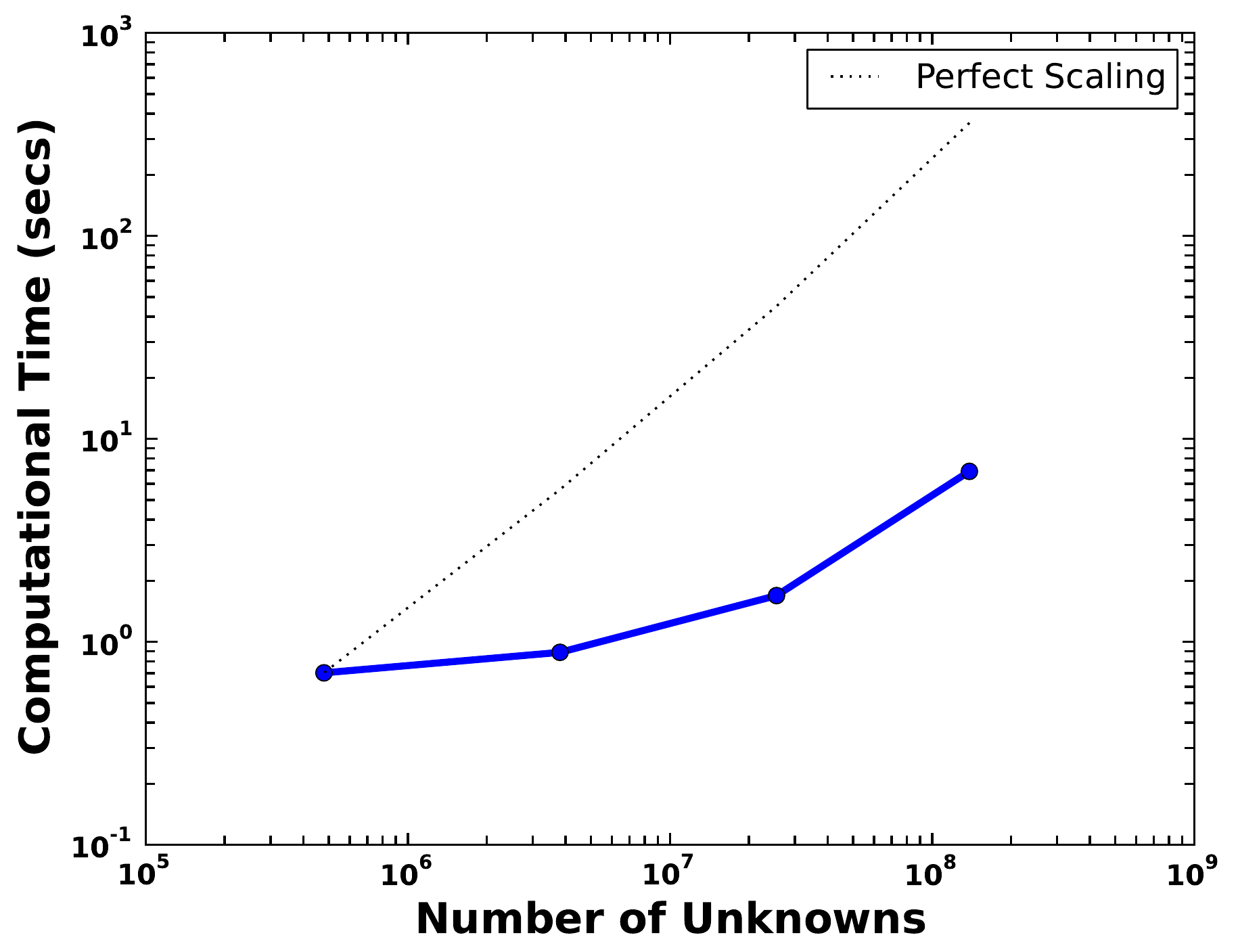}
   \caption{Average sample time}
   \label{subfig:SPE10_time}
  \end{subfigure}
  \begin{subfigure}{0.49\textwidth}
    \centering
    \includegraphics[scale=.33]{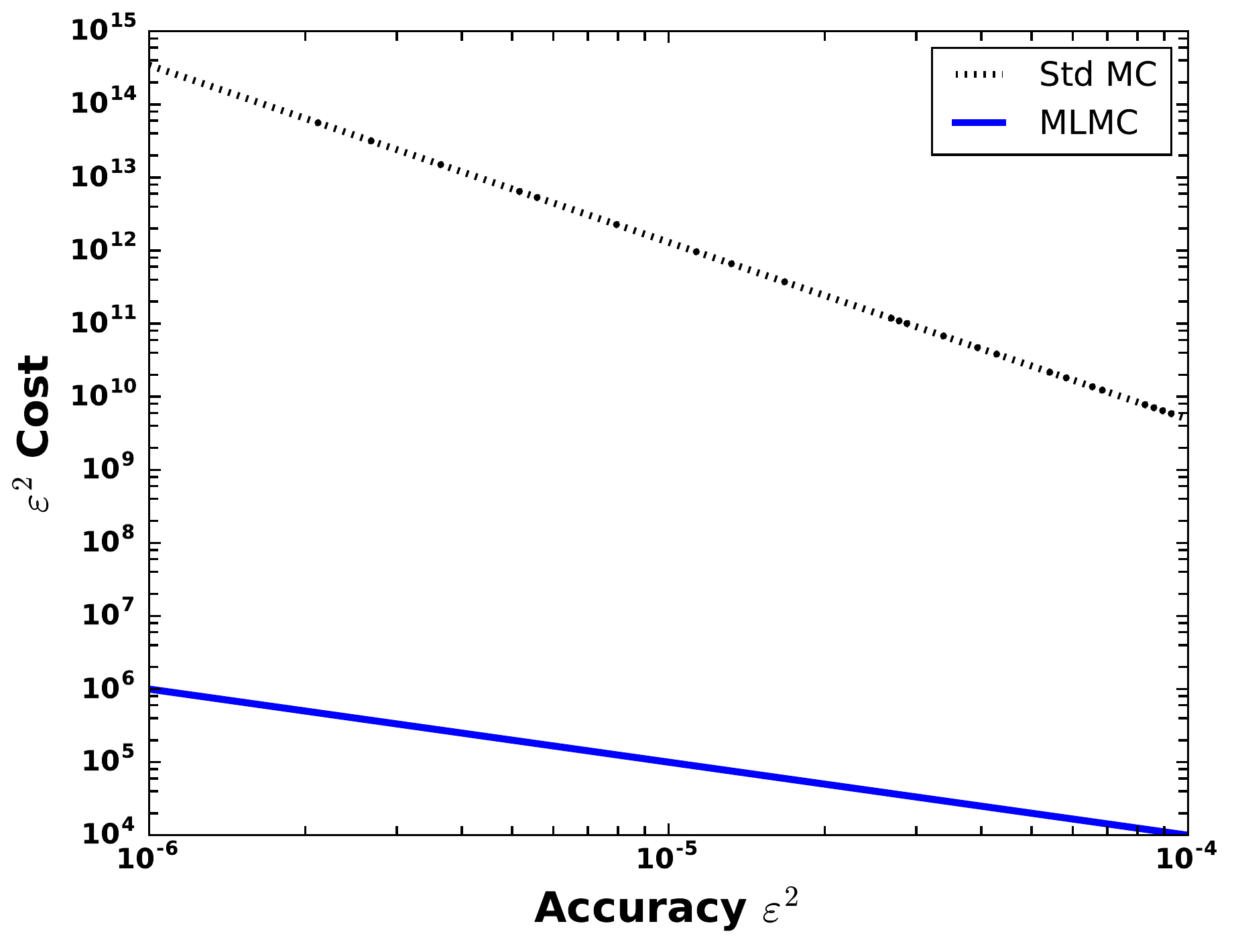}
   \caption{Predicted Asymptotic Order of Cost}
   \label{subfig:SPE10_theoretical}
  \end{subfigure}
  \begin{subfigure}{0.99\textwidth}
  \begin{center}
  \scalebox{0.9}{
  \begin{tabular}{lrrrr}
  \toprule
  Level &  0 (fine)&  1&  2& 3 (coarse)  \\
  \midrule
  Number of samples & 148 & 419 & 766 & 4097 \\
  \hline 
  Realize input coefficients (SPDE sampler) & 230.8140 &110.5328 & 74.4681 & 143.3106 \\
  Assemble systems \eqref{eq:ml-mixedDarcyDiscrete} on levels $\ell$ and $\ell+1$ & 144.0835 &28.6440 &7.4916 & 8.4634  \\
  Build AMG preconditioner for $\tilde{S}_{\ell}$ and $\tilde{S}_{\ell+1}$ in \eqref{eq:ldu} & 2.9663 & 1.1234 & 0.6655 & 1.3899 \\
  Solve \eqref{eq:ml-mixedDarcyDiscrete} on levels $\ell$ and $\ell+1$& 960.1115 & 1000.1379 & 1152.0609 & 2845.5386\\
  \hline
   Total time & 1337.9754 &  1140.4381 & 1234.6861 & 2998.7025\\
  \bottomrule
  \end{tabular} }
  \end{center}
    \caption{
    Detailed computational times (secs) for each level of the MLMC simulation. 
  }
  \label{subfig:SPE10_mlmc_table}
  \end{subfigure}
  \caption{MLMC results for SPE10 problem where the target MSE is
      $\varepsilon^2=6.25\cdot 10^{-6}$ using $2304$ processes with
      approximately $1.4 \cdot 10^8$ velocity and pressure degrees of freedom
      on the fine level. The MC estimator at each level
      and the MC estimator of the correction is shown in (a) for each level, and the 
      variance reduction of the multilevel method is demonstrated in (b),
      where the variance of the MC estimator and MC estimator of the correction 
      term are plotted for each level.
      Plot (c) shows the average sampling
      time to generate the required Gaussian field
  realizations and solve the forward model for each level versus the number of
  unknowns. The theoretical asymptotic order of cost to
    achieve a MSE of $\varepsilon^2$ for this problem formulation comparing the
standard
  MC and MLMC estimators is shown in Plot (d), where it is visible that the MLMC
  method leads to significant computational savings over the standard MC
  estimator. The time spend computing on each level is shown in (e) where 
  the majority of time is spent on the coarse grid.}\label{fig:SPE10_mlmc}
\end{figure}

We now consider incorporating data from the SPE10 benchmark into a MLMC simulation.
The random permeability coefficient $k(\bx,\omega)$ is modeled as
log-normal random field with mean equal to the absolute permeability given by
the SPE10 dataset, so that
$$\exp[\log[k_{SPE10}(\bx)] + \theta(\bx,\omega)],$$ where $\theta(\bx,\omega)$
has an exponential covariance with $\sigma^2=1$ and correlation length
$\lambda=50$(ft). The hierarchical SPDE sampler with non-matching mesh embedding is
used to generate the realizations of the random field $\theta(\bx,\omega)$.
Figure~\ref{fig:spe10_perm} illustrates the permeability field of the SPE10
dataset $k_{SPE10}(\bx)$.
The forward problem is given by 
~\eqref{eq:darcy} with boundary conditions
  \begin{equation*}
    \begin{cases}
      -p&=1 \quad \text{on} \ \Gamma_{in}=\{0\}\times (0,2200) \times (0,170), \\
      -p&=0 \quad \text{on} \ \Gamma_{out}=\{1200\}\times (0,2200) \times (0,170), \\
      \bq \cdot {\bf n} &= 0\quad  \text{on} \ \Gamma_{s} := \partial D \setminus \left(\Gamma_{in} \cup \Gamma_{out}\right). \\
    \end{cases}
  \end{equation*}
    The quantity of interest is the expected value of the Darcy pressure
    evaluated on the fine element around the point $\bx^*=(600,1100,85).$
Next we demonstrate the performance of a $4$-level MLMC simulation using $2304$ processes with
approximately $1.4 \cdot 10^8$ velocity and pressure degrees of freedom on the fine level in
Figure~\ref{fig:SPE10_mlmc}. The target MSE is set to
$\varepsilon^2=6.25 \cdot 10^{-6}$.

The MC estimator of the QOI and the correction term on each level 
and the variance of both quantities are displayed
in Figures~\ref{subfig:SPE10_mean} and \ref{subfig:SPE10_var}, respectively.
The blue solid line represents the standard MC estimator at each level $\ell$, whereas
the green dashed line represents the MC estimator of the correction term.
The MC estimates demonstrate the expected behavior, 
confirming the benefits of the multilevel approach, where the variance 
of the correction term decreases as the spatial resolution increases.
The average sampling time to generate the required Gaussian field realizations
and solve the forward model at each level is shown in
Figure~\ref{subfig:SPE10_time}.
This plot demonstrates the desired scalability of
the solution strategy for the forward problem and of the proposed hierarchical sampler with non-matching
mesh embedding.
The predicted $\varepsilon^2$-cost of
the standard MC and MLMC estimators is shown in
Figure~\ref{subfig:SPE10_theoretical} using Theorem 2.3
from~\cite{teckentrup2013further} with
numerically observed constants estimated from
Figures~\ref{subfig:SPE10_mean}-\ref{subfig:SPE10_time}.
In Figure~\ref{subfig:SPE10_mlmc_table}, the computational time spent on 
each level illustrates the expected performance, where the majority of
time is spent generating samples on the coarsest level.
Next we consider the scalability of adaptive MLMC simulations for the SPE10 problem.
Numerical results are presented in Table~\ref{tab:spe10_weak_table} for
 three different MLMC simulations where the spatial resolution
is increasing, while the number of degrees of freedom per process remains approximately fixed.
The desired MSE is chosen to balance the estimated
discretization error with the sampling error, resulting in more necessary samples
for a finer spatial discretization while the
majority of the samples are generated on the coarsest level.

\begin{table}[tbhp]
  \caption{Weak scaling of SPE10 problem for adaptive MLMC simulations with 
  approximately $6.2\cdot10^4$ pressure/velocity degrees of freedom
  per process on the fine level.
  In each row, the number of MPI processes are listed along with the
  number of velocity and pressure degrees of freedom on the
  fine level ($\ell=0$), the target MSE $\varepsilon^2$, 
  the total computational wall time to run the MLMC simulation,
  the number of samples computed on the fine level ($N_0$), and the total
  number of computed samples on all levels for a 4 levels method.}
  \label{tab:spe10_weak_table}
  \begin{center}
  \begin{tabular}{lllllr}
    \toprule
    Processes&  DOF ($\ell=0$)& Target $\varepsilon^2$ & Wall Time (s)&  $N_0$&  Total Samples \\
    \midrule
    36 &   2.19$\cdot 10^6$& $1.00\cdot 10^{-4}$ & 3.09$\cdot10^2$ &  12&   641 \\
    288 &  1.74$\cdot 10^7$& $2.50\cdot 10^{-5}$ & 1.26$\cdot10^3$ &  46&  2157 \\
    2304 & 1.39$\cdot 10^8$& $6.25\cdot 10^{-6}$ & 6.45$\cdot10^3$ &  148&   5430\\
    \bottomrule
  \end{tabular}
  \end{center}

\end{table}

\subsection{Discussion}\label{subsec:discussion}
These results demonstrate that the parallel hierarchical SPDE sampler
with non-matching domain embedding coupled with a scalable forward model solver
allows for accurate large-scale MLMC simulations to be performed, 
which otherwise would not have been feasible.
We have primarily focused on a novel, hierarchical sampling 
technique with non-matching mesh embedded for generating the necessary 
realizations of a Gaussian random fields. Our implementation is highly  
scalable with respect to the number of degree of freedom on the fine grid (74\%
parallel efficiency on 10 thousand processors). In addition, the overall parallel
efficiency of the MLMC methods can be further improved by exploiting an additional 
layer of parallelism to generate multiple independent samples concurrently.
In the current implementation, each sample is computed sequentially.
We only exploit parallelism in the spatial dimension and, therefore, we use
the same number of processors to solve the fine and coarser problems.
This leads to an under utilization of computational resources on the coarser levels,
where the problem size is too small with respect to the number of processors 
employed, and causes deterioration of performance in the linear solver phase. 
A significant improvement in the scalability of our methodology requires 
repartitioning the coarser problems on a subset of processes, so that multiple coarse samples
can be computed in parallel and asynchronously by different subsets of processes.
This approach has been investigated in \cite{GmeinerDRSW16}, and can be applied directly applied 
also to our methods, modulo some implementation challenges and nuances all left for possible future studies.

\section{Conclusions} \label{sec:conclusions} 
The ability to efficiently generate samples of a Gaussian random field at different
spatial resolutions is an essential component of large-scale sampling-based methods for forward
propagation of uncertainty. 
We propose a hierarchical sampling method based on the solution of a reaction-diffusion stochastic
PDE using a domain embedding technique with two non-matching meshes. 
The stochastic PDE is discretized and solved on a regular domain with a
structured mesh, then transferred to the original, unstructured mesh of
interest.
The proposed sampling method allows for Gaussian random field realizations
to be scalably generated for complex spatial domains, by leveraging efficient 
preconditioning techniques for the iterative solution
of the discrete saddle-point problem arising from the discretization of 
the stochastic PDE on a regular domain.  
A hierarchical version of this process is explored, and numerical
results demonstrate the scalability of the proposed method for generating
realizations of a log-normal random field for large-scale simulations of flow in porous media. 
Additionally the sampling method is used in MLMC simulations
of subsurface flow problems and numerical 
results are presented, which demonstrate the scalability of the hierarchical SPDE
sampler
with non-matching domain embedding for MLMC simulations of subsurface flow
problems with over 470 million parameters in the stochastic dimension and 
1.9 billion spatial unknowns.

\section*{Acknowledgments}
R.K. and P.Z. acknowledge the support by the Swiss Commission for Technology and
Innovation via the SCCER-FURIES and by the Swiss National Science Foundation,
via the projects ExaSolvers - Extreme Scale Solvers for Coupled Systems, and
``Geometry-Aware FEM in Computational Mechanics''.

\end{document}